\documentclass[12pt]{amsart}
\title{Edited $4\Theta$-embeddings of
Jacobians}
\author{Greg W. Anderson}
\address{University of Minnesota, MN 55455 USA}
\email{gwanders@math.umn.edu}
\date{This paper is to appear in the Michigan Mathematical Journal.
It was accepted pending minor revisions on Mar 21, 2003. 
Those revisions are now complete. 
This is the final draft submitted for publication.}

\thanks{2000 {\it Mathematics Subject Classification}. Primary: 14H40, 14H42}
\newcommand{\TT}{{\mathcal T}}
\newcommand{\HH}{{\mathcal H}}
\newcommand{\OO}{{\mathcal O}}
\newcommand{\EE}{{\mathcal E}}
\DeclareMathOperator{\supp}{{\mathrm{supp}}}
\DeclareMathOperator{\abel}{{\mathrm{abel}}}
\newcommand{\card}{{\#}}
\newtheorem{Lemma}[subsubsection]{Lemma}
\newtheorem{Theorem}[subsubsection]{Theorem}
\newtheorem{Corollary}[subsubsection]{Corollary}
\newtheorem{Proposition}[subsubsection]{Proposition}
\newcommand{\ZZ}{{\mathbb Z}}
\newcommand{\CC}{{\mathbb C}}
\newcommand{\RR}{{\mathbb R}}

\newcommand{\barpartial}{{\bar{\partial}}}
\DeclareMathOperator{\Aut}{{\mathrm{Aut}}}
\begin{document}
\maketitle
\begin{abstract} By the Lefschetz embedding theorem
a principally polarized $g$-dimensional abelian variety is embedded into
projective space by the linear system of
$4^g$ half-characteristic theta functions. Suppose we {\em edit}
this linear system by dropping all the theta functions vanishing at the
origin to order greater than parity requires. We prove that for Jacobians
the edited $4\Theta$ linear system still defines an embedding into
projective space. Moreover, we prove that the
projective models of Jacobians arising from the elementary algebraic  construction of
Jacobians recently given by the author are (after passage to linear hulls)
copies of the edited
$4\Theta$ model. We obtain our results by aptly combining the
quartic and determinantal identities satisfied by the Riemann theta
function. We take the
somewhat nonstandard tack of working in the framework of Weil's old book on
K\"{a}hler varieties in order to avoid having to make
extremely complicated calculations.
\end{abstract}

\section{Introduction}

The point of departure for this paper is the elementary
algebraic construction of Jacobians given in
\cite{Anderson}. We begin by reviewing that
construction. For brevity's sake we reformulate the construction in
terms of line bundles rather than divisors. Let
$X$ be a nonsingular projective algebraic curve of genus $g>0$.
Although the main work of this paper takes place over the complex
numbers, for the moment we take as ground field any algebraically
closed field.  Fix an
integer
$n\geq g+2$. For $i=0,\dots,n+1$, let $f\mapsto f^{(i)}$ denote the
operation of pull-back via the $i^{th}$ projection
$X^{\{0,\dots,n+1\}}\rightarrow X$. Fix a line bundle $\EE$ on $X$ of
degree $n+g-1$.  Given any line bundle
$\TT$ on
$X$ of degree
$0$, let
$u$ (resp. $v$) be a row vector of length $n$ with entries
forming a basis of $H^0(X,\TT^{-1}\otimes \EE)$ 
(resp., $H^0(X,\TT\otimes \EE)$) over the ground field,
and let $\abel(\TT)$ be the $n$ by $n$ matrix with entries
$$\abel(\TT)_{ij}:=\left|\begin{array}{c}
\widehat{v^{(0)}}\\
\vdots\\
\widehat{v^{(i)}}\\
\vdots\\
\end{array}\right| \cdot
\left|\begin{array}{c}
\vdots\\\widehat{u^{(i)}}\\\vdots\\\widehat{u^{(n+1)}}
\end{array}\right|\cdot
\left|\begin{array}{c}\vdots\\
\widehat{v^{(j)}}\\
\vdots\\
\widehat{v^{(n+1)}}
\end{array}\right|\cdot
\left|\begin{array}{c}
\widehat{u^{(0)}}\\\vdots\\\widehat{u^{(j)}}\\\vdots
\end{array}\right|,
$$
where the leftmost determinant is that obtained
by stacking the row vectors $v^{(i)}$ to form an $n+2$ by $n$ matrix
with rows numbered from $0$ to $n+1$, 
then striking row $0$ and row $i$ to obtain a square matrix,
and finally
taking the determinant; the other determinants are analogously formed.
Up to a nonzero scalar multiple the matrix $\abel(\TT)$ is  
independent of the choice of bases
$u$ and
$v$ and moreover depends only on the isomorphism class of the line bundle
$\TT$. It is easy to see that
$\abel(\TT)$ does not vanish identically.   The construction 
$\TT\mapsto \abel(\TT)$ maps the set of isomorphism classes of degree
zero line bundles on $X$ to the projective space of lines in
the space of $n$ by $n$ matrices with entry in $i^{th}$ row
and $j^{th}$ column drawn from the space
$$
H^0\left(X^{\{0,\dots,n+1\}},\frac{\displaystyle\bigotimes_{\ell=0}^{n+1}
\left(\EE^{(\ell)}\right)^{\otimes 4}}
{\left(\EE^{(0)}\right)^{\otimes 2}
\otimes\left(\EE^{(i)}\right)^{\otimes
2}
\otimes\left(\EE^{(j)}\right)^{\otimes
2}
\otimes
\left(\EE^{(n+1)}\right)^{\otimes 2}}\right).
$$
In
\cite{Anderson}:
\begin{itemize}
\item the map $\TT\mapsto \abel(\TT)$ is shown to be injective,
\item 
an
explicit set of homogeneous equations cutting out the image
of the map
$\TT\mapsto \abel(\TT)$ is exhibited, and
\item the image  of the map $\TT\mapsto \abel(\TT)$ is
explicitly equipped with an algebraic group law commuting with the
tensor product of degree zero line bundles. 
\end{itemize}
Thus the Jacobian of $X$ is
constructed in elementary algebraic fashion.

The question motivating this paper
is the following:
\begin{quote}
What linear system of effective divisors of the Jacobian
arises from the projective embedding $\TT\mapsto \abel(\TT)$?
\end{quote}
We attack the question by the methods of complex algebraic geometry
and obtain a complete answer in that setting.
The question remains open in
positive characteristic.

The answer we finally obtain is surprisingly simple. In the
classical-style language of \cite{Fay},
\cite{MumfordTataI} and
\cite{MumfordTataII} the answer takes the following form. The curve $X$
is now a compact Riemann surface. Let
$\tau$ be the period matrix of $X$. Let $M$ be the set of column vectors
of length
$2g$ with entries in the set $\{0,1/2\}$. We write such column vectors in
block form
$\left[\begin{array}{c}a\\b\end{array}\right]$
where $a$ and $b$ are both of length $g$. Let $M_{\leq 1}$
be the subset of $M$ consisting of those
$\left[\begin{array}{c}a\\b\end{array}\right]$ such that the corresponding
half-characteristic classical theta function
$$\theta\left[\begin{array}{c}a\\b\end{array}\right]
(w,\tau ):=
\sum_{\ell\in \ZZ^g}
\exp\left(\pi i(\ell+a)^T\tau (\ell+a)+2\pi
i(\ell+a)^T(w+b)\right)$$
vanishes at $w=0$ to order not exceeding $1$, i.~e., to order not greater
than that dictated by parity considerations. Then---independent of the
choice of
$n$ and the line bundle
$\EE$---it turns out that the linear system we are looking for is the {\em
edited}
$4\Theta$ linear system
$$\left\{\left.\theta\left[\begin{array}{c}a\\b\end{array}\right]
(2w,\tau
)\right|\left[\begin{array}{c}a\\b\end{array}\right]\in
M_{\leq 1}\right\}.$$
(Here and throughout the paper we abuse language by identifying linear
systems of effective divisors with linearly independent sets
of sections of line bundles in the obvious way.) 
Of course, for this answer to make sense, the edited
$4\Theta$ linear system has to embed the Jacobian into projective
space, i.~e., the ``edited version'' of the Lefschetz
embedding theorem has to hold. The latter we prove in this
paper (see Theorem~\ref{Theorem:EditedFourTheta} below) by complex
analytic methods independent of---but largely parallel to---the methods of
\cite{Anderson}. 
The theorem has some
content since, for example, $M_{\leq 1}\neq M$ for all
hyperelliptic curves of large genus.  
We conclude the paper by explaining in detail how to factor
the map $\TT\mapsto \abel(\TT)$ through the edited $4\Theta$-embedding,
thus fully answering the motivating question over the complex numbers.

The main technical problem faced in this paper is that of aptly expressing
and combining the quartic and determinantal identities satisfied by the
Riemann theta function. To handle the problem we take the somewhat
nonstandard tack of working in the framework of Weil's old book
\cite{Weil}  on K\"{a}hler varieties. 
This is advantageous because Weil's austere conceptual
approach to theta functions obviates quite a lot of bookkeeping.
For example, there is no need to choose $A$- and $B$-cycles,
and hence there is a generally lower risk of making a 
sign error, cf.\ \cite[p.\ 3.81]{MumfordTataII}.
 In \S\ref{section:ThetaTools} below we explain our
Weil-style point of view on classical theta identities.   In
\S\ref{section:AbeliantIdentity} below we fit the quartic and
determinantal theta identities together within the Weil framework,
thereby obtaining our main results.

\section{A theta function toolkit}\label{section:ThetaTools}  We
review what we need of the general theory of theta functions, 
following 
\cite{Weil}. Then we work through a series of examples in order to
bring notions treated in
\cite{Fay}, \cite{GriffithsHarris}, \cite{MumfordTataI},
\cite{MumfordTataII} and
\cite{MumfordTataIII} into the Weil picture.

\subsection{A general theta formalism}

\subsubsection{The setting}\label{subsubsection:Setting}
Fix a  compact complex manifold $V$ of K\"{a}hler type.
(In applications $V$ is going to be an abelian
variety, or a compact Riemann surface, or a product of such.)  Fix a
universal covering map
$\tilde{V}\rightarrow V$ (we do not bother to give the map a name) and
denote its automorphism group by
$G$. The group law in $G$ is by definition composition of functions
and hence $G$ acts naturally on the left of $\tilde{V}$.  Given
$\tilde{v}\in
\tilde{V}$ and its image $v\in V$ under
the covering map
$\tilde{V}\rightarrow V$, we call $\tilde{v}$ a {\em lifting}
of $v$. Analogously
we speak of liftings of paths. Given a differential
form
$\omega$ on
$V$ and its pull-back
$\tilde{\omega}$ via the covering map
$\tilde{V}\rightarrow V$, we call $\tilde{\omega}$
the {\em lifting} of $\omega$.  Given a closed
$1$-form
$\zeta$ on $V$ and a function
$z$ on $\tilde{V}$ such that 
$dz$ is the lifting of $\zeta$, we call
$z$ a {\em primitive} of
$\zeta$. 

\subsubsection{Systems of multipliers} We call a family
$\{F_\sigma\}$ of \linebreak nowhere-vanishing
holomorphic functions on $\tilde{V}$ indexed by
$\sigma\in G$
a {\em system of multipliers} under the
following conditions:
\begin{itemize}
\item $F_{\sigma\tau}=\tau^*F_\sigma\cdot F_\tau$
for all $\sigma,\tau\in G$.
\item $d\log F_\sigma$ is the lifting of a holomorphic
$1$-form on $V$ for all $\sigma\in G$.
\end{itemize}
When the roles of $\tilde{V}$ and $G$ require
emphasis, we say that $\{F_\sigma\}$ is a system on
multipliers {\em on $\tilde{V}$ relative to $G$}. If all the functions
$F_\sigma$ are constants, we say that $\{F_\sigma\}$ is a system of {\em
constant} multipliers. A system of 
constant multipliers is simply a
homomorphism from $G$ to the group of nonzero complex
numbers.  If all the functions
$|F_\sigma|$  are identically equal to $1$, we say that
$\{F_\sigma\}$ is {\em unitary}.
By the
Maximum Principle a unitary
system of multipliers is a system of
constant multipliers.
    We say that systems of multipliers
$\{F_\sigma\}$ and
$\{F'_\sigma\}$ are {\em equivalent} if for some
nowhere-vanishing holomorphic function $u$ on
$\tilde{V}$ we have
$F'_\sigma/F_\sigma =\sigma^*u/u$ for
all
$\sigma\in G$. 
By Hodge theory every system
of constant multipliers is equivalent to a unitary
system of multipliers. Moreover, distinct unitary systems of multipliers
are inequivalent.

\subsubsection{Theta functions}
\label{subsubsection:ThetaDef}
A not-identically-vanishing
meromorphic function $\vartheta$ on
$\tilde{V}$ transforming for some 
 system of multipliers
$\{F_\sigma\}$ by the rule
$\sigma^*\vartheta=F_\sigma\cdot \vartheta$
for all $\sigma\in G$ is called a {\em theta function}.
When the roles of $\tilde{V}$ and $G$ require
emphasis, we say that $\vartheta$ is a theta function {\em on
$\tilde{V}$ relative to $G$}.  Given a theta function $\vartheta$ 
transforming according to a system of multiplier $\{F_\sigma\}$, we
say that
$\vartheta$ {\em determines} the system of multipliers
$\{F_\sigma\}$; clearly $\{F_\sigma\}$ is uniquely determined by
$\vartheta$. The divisor of a theta function
$\vartheta$ is
$G$-invariant and hence descends to a divisor of
$V$, say $D$; in this situation we say that
$\vartheta$ {\em represents} $D$,
that $\vartheta$ is {\em effective}
if $D$ is effective, and that $\vartheta$ is
{\em trivial} if $D=0$.  We say
that theta functions are {\em equivalent} if they
represent the same divisor of $V$. Equivalent
theta functions determine equivalent multiplier
systems. By the Maximum Principle an effective theta
function with
$G$-invariant absolute value is constant.

\subsubsection{The
first
Chern class of a
divisor}\label{subsubsection:DivisorCohomologyClass}
 Fix a divisor
$D$ of $V$. For the purpose of checking signs and factors of $2$
and $\pi$, we briefly recall the method used in
\cite[Chap.\ 5]{Weil} for obtaining the first Chern class $c_1(D)$ in the
de Rham cohomology of $V$. Fix an open
covering
$\{U_i\}$ of $V$ by nonempty open sets such that
that
$D$ restricted to $U_i$ is the divisor of a meromorphic
function
$f_i$ on $U_i$.  For indices $i$ and $j$ such that
$U_i\cap U_j\neq
\emptyset$, let
$F_{ij}$ be the unique nowhere-vanishing holomorphic
function on
$U_i\cap U_j$ such that 
$$f_j\vert_{U_i\cap U_j}=F_{ij}\cdot f_i\vert_{U_i\cap
U_j}.$$
The
family
$\{F_{ij}\}$ thus defined is called the {\em system of 
transition functions} associated to
$\{U_i\}$ and $\{f_i\}$.  By a method recalled in the course
of the proof of Proposition~\ref{Proposition:PoincareRecognition} below,
it is possible to construct for each $i$ a smooth
$(1,0)$-form
$\eta_i$ on $U_i$  such that  
$$\frac{1}{2\pi i}d\log
F_{ij}=\eta_j\vert_{U_i\cap U_j}-\,
\eta_i\vert_{U_i\cap U_j}$$
whenever $U_i\cap U_j\neq \emptyset$;
any such family
$\{\eta_i\}$ will be called a {\em connection} for 
$\{F_{ij}\}$. Given any connection $\{\eta_i\}$ for $\{F_{ij}\}$, there
exists a unique smooth closed
$2$-form
$\alpha$ on $V$ described locally by the conditions
$$\alpha\vert_{U_i}=d\eta_i.$$ 
The closed
$2$-form
$\alpha$ is called the {\em curvature} of the connection
$\{\eta_i\}$. The  de Rham cohomology class of the $2$-form
$\alpha$ depends only on
$D$, not on the intervening choices, and this class is none other than
$c_1(D)$.

\begin{Lemma}\label{Lemma:LocalProblem}
Let $\phi$ be a smooth compactly supported function on
$\CC^n$.  Let $\beta$ be a smooth closed $(2n-2)$-form
defined on an open set $U\subset\CC^n$ containing
the support of $\phi$.  The integral
$\int_{U\setminus \{z_1=0\}} d\phi\wedge
d\log z_1
\wedge \beta$
is absolutely convergent and equals
$-2\pi i
\int_{U\cap\{z_1=0\}}\phi\beta$.
\end{Lemma}
\proof Sophomore calculus.
\qed

\begin{Proposition}
\label{Proposition:PoincareRecognition}
Let $n$ be the complex dimension of $V$. Let $D$ be a
complex submanifold of
$V$  of codimension $1$. A closed
$2$-form $\alpha$ on $V$ belongs to the cohomology
class
$c_1(D)$ if and only if
$\int_V\alpha\wedge
\beta=\int_D
\beta$ for all closed $(2n-2)$-forms $\beta$ on $V$.
\end{Proposition}
\proof By Poincar\'{e} duality it suffices merely to
exhibit a closed
$2$-form
$\alpha$ belonging to the class $c_1(D)$ such that $\int_V\alpha\wedge
\beta=\int_D
\beta$ for all closed $(2n-2)$-forms $\beta$ on $V$, and it is well
known in principle how to do this. We take care over the details just
for the purpose of checking signs and factors of $2$ and $\pi$. By
hypothesis we can choose a finite open covering
$\{U_i\}$ of $V$ by coordinate patches
and for each $i$ a holomorphic function $f_i$ on $U_i$ 
belonging to some coordinate system on $U_i$
such that $D\cap U_i= \{f_i=0\}$. Let $\{F_{ij}\}$
be the system of transition functions associated to 
$\{U_i\}$ and $\{f_i\}$. Fix a partition of unity
$\{\phi_i\}$ subordinate to $\{U_i\}$. 
For all indices $i$ and $\ell$  there
exists a unique smooth
$(1,0)$-form
$\eta_{i\ell}$ on
$U_i$ 
such that 
$$
2\pi
i\cdot\eta_{i\ell}\vert_{U_i\cap
U_\ell}=\phi_\ell\vert_{U_i\cap U_\ell}\cdot d \log
F_{\ell i}\;\;\mbox{if $U_i\cap U_\ell\neq
\emptyset$},\;\;\;\;\eta_{i\ell}\vert_{U_i\setminus\supp
\phi_\ell}=0.$$ 
For all
indices
$\ell$ there exists a unique smooth $(1,0)$-form
$\zeta_\ell$ on \linebreak $V\setminus \supp D$ such
that
$$\zeta_\ell\vert_{U_\ell\setminus \supp D}=
\phi_\ell\vert_{U_\ell\setminus \supp D}
\cdot d\log f_\ell\vert_{U_\ell\setminus
\supp D},\;\;\;
\zeta_{\ell}\vert_{V\setminus (\supp D\cup \supp
\phi_\ell)}=0.$$ 
The family $\left\{\sum_{\ell}
\eta_{i\ell}\right\}$ is then connection for $\{F_{ij}\}$
with curvature
$\alpha$ satisfying the relation
$$-2\pi i\cdot\alpha\vert_{V\setminus
D} =\sum_\ell d\zeta_\ell.
$$
The closed $2$-form $\alpha$ has the desired
property by Lemma~\ref{Lemma:LocalProblem}.
\qed

\subsubsection{Weil gauges}
  We say that a
real-valued function $\Phi$ on $\tilde{V}$ is a {\em
Weil gauge}  if of the form
$\Phi=\sum_{i=1}^N c_i|z_i|^2$
where the $z_i$ are
primitives  of holomorphic $1$-forms
on $V$ and the $c_i$ are real constants. The
notion of Weil gauge plays a key if
anonymous role in 
\cite[Chap.\ V]{Weil}.

\begin{Theorem}\label{Theorem:ThetaTool}
Fix a divisor
$D$  of $V$ and a Weil gauge $\Phi$ on $\tilde{V}$.
The following conditions are equivalent:
\begin{enumerate}
\item $\frac{i}{2}\partial\barpartial\Phi$ is
the lifting of a closed real
$(1,1)$-form on $V$
belonging to the first Chern class $c_1(D)$.
\item There exists a theta function 
$\vartheta$ unique up to a nonzero constant factor 
such that $\vartheta$ represents $D$ 
and 
$e^{-\pi\Phi}|\vartheta|^2$ is 
$G$-invariant. 
\end{enumerate}
(If the first condition holds we say that
$\Phi$ is a {\em gauge} for $D$;
Proposition~\ref{Proposition:PoincareRecognition}  is sometimes
convenient for checking this. Under the condition that
$e^{-\pi\Phi}|\vartheta|^2$ is
$G$-invariant we say that
$\vartheta$ is {\em $\Phi$-normalized}.
The ratio of any two $\Phi$-normalized theta functions
necessarily transforms according to a unitary character
of
$G$.)
\end{Theorem}
\subsubsection{Partial sketch of proof} The
Maximum Principle proves uniqueness. The proofs of the
implications (1$\Rightarrow$2) and (2$\Rightarrow$1)
are more or less evident modifications of the
proofs of
\cite[Thm.~2, Chap.~V]{Weil} and
\cite[Prop.~3, Chap.~V]{Weil}, respectively. We therefore need not provide
a detailed proof of the equivalence  (1$\Leftrightarrow$2).
Still, since we have
superficially modified Weil's theory  by stressing the notion of
Weil gauge, we are under some obligation at least to check signs and factors
of
$2$ and
$\pi$. We therefore compromise: we very rapidly sketch a
proof of the implication (2$\Rightarrow$1) and then omit further
details.

 Let $\{F_\sigma\}$ be the system of multipliers determined by
$\vartheta$. 
 By
hypothesis we have
$$|F_\sigma|^2=e^{\pi(\sigma^*\Phi-\Phi)}$$
 and hence
$$d\log F_\sigma=\pi
\partial(\sigma^*\Phi-\Phi)$$
for all $\sigma\in G$.
Let
$\{U_i\}$ be a finite covering of $V$ by geodesically convex nonempty
open sets. For each $i$ fix a section
$s_i:U_i\rightarrow
\tilde{V}$ of the covering map. Let $\{F_{ij}\}$ be the system
of transition functions associated to $\{U_i\}$ and
$\{s_i^*\vartheta\}$---meaning if $U_i\cap U_j\neq
\emptyset$, then we have
$$s_j^*\vartheta\mid_{U_i\cap U_j}=F_{ij}\cdot
s_i^*\vartheta\mid_{U_i\cap U_j}.$$
 Whenever
$U_i\cap U_j\neq
\emptyset$, choose $\sigma_{ij}\in G$ to satisfy the
condition
$$
s_j\vert_{U_i\cap U_j}=\sigma_{ij}\circ s_i\vert_{U_i\cap U_j}.$$
Such a choice exists and is unique because $U_i\cap U_j$
is geodesically convex and {\em a fortiori} connected.
Then whenever $U_i\cap U_j\neq \emptyset$ we have
$$\begin{array}{rcl}
\frac{1}{2\pi i}d\log F_{ij}
&=&\frac{1}{2\pi i}s_i^*d\log F_{\sigma_{ij}}\vert_{U_i\cap U_j}\\\\
&=&\frac{1}{2i}s_i^*\partial
(\sigma_{ij}^*\Phi-\Phi)\vert_{U_i\cap U_j}\\\\
&=&\frac{1}{2i}s_j^*\partial\Phi\mid_{U_i\cap U_j}
-\frac{1}{2i}s_i^*\partial\Phi\mid_{U_i\cap U_j},
\end{array}$$
i.\ e., $\left\{\frac{1}{2i}s_i^*\partial \Phi\right\}$
is a connection for $\{F_{ij}\}$. Therefore the real closed
$(1,1)$-form
$\alpha$ defined locally by the conditions 
$$\alpha\mid_{U_i}=
\frac{i}{2}s_i^*\partial\barpartial\Phi$$ belongs to the class
$c_1(D)$.
\qed

\subsubsection{Complement}
\label{subsubsection:BundleGauge}
 With an emphasis on line
bundles rather than on divisors, Theorem~\ref{Theorem:ThetaTool} takes
the following form. Let
$\EE$ be a line bundle on
$V$ with pullback
$\tilde{\EE}$ to $\tilde{V}$. Let $\Phi$ be a Weil gauge on
$\tilde{V}$. Then the following conditions are equivalent:
\begin{enumerate}
\item 
$-\frac{i}{2}\partial\bar{\partial}
\Phi$ is the lifting of a closed real $(1,1)$-form on $V$
belonging to the first Chern class of the line bundle $\EE$.
\item There exists a multiplier system $\{F_\sigma\}$
such that for some global trivialization $\tilde{e}$
of $\tilde{\EE}$ we have
$$\sigma^*\tilde{e}=F_\sigma \cdot \tilde{e},\;\;\; 
|F_\sigma|^2=e^{\pi(\sigma^*\Phi-\Phi)}$$ for all
$\sigma\in G$, and moreover the multiplier system
$\{F_\sigma\}$ thus attached to $\EE$ and $\Phi$ is unique.
\end{enumerate}
Further, with $\EE$, $\Phi$ and
$\{F_\sigma\}$ as above, if
$\EE=\OO_V(-D)$ for some divisor $D$, then $\Phi$ is a gauge for $D$ and
$\{F_\sigma\}$ is the multiplier system determined by any
$\Phi$-normalized theta function representing $D$.

\subsection{Example: Principally polarized complex tori}
\label{subsection:PrincipallyPolarized}

\subsubsection{Definition}
Following
\cite[Chap.\ VI]{Weil}, we define a {\em principally polarized complex
torus} of complex dimension
$g$ to be a triple 
$$(W,H,\Lambda)$$
 consisting of
\begin{itemize}
\item a $g$-dimensional complex vector space $W$,
\item a  positive definite hermitian form 
$$H:W\times W\rightarrow\CC$$
anti-linear on the left  and
\item a cocompact discrete subgroup $\Lambda\subset W$
\end{itemize}
such that
\begin{itemize}
\item the relations
$$\Im
H(A_i,B_j)=\delta_{ij},\;\;\;
\Im H(A_i,A_j)=0=\Im H(B_i,B_j)
$$
hold for at least one
$\ZZ$-basis  $\{A_i,B_i\}_{i=1}^g$ of $\Lambda$.
\end{itemize}
Any $\ZZ$-basis for $\Lambda$ with the special property above we call
{\em symplectic}. The hermitian form
$H$ naturally gives rise to an invariant
K\"{a}hler metric on the complex Lie group $W/\Lambda$.

\subsubsection{Specialization of the setting}
Fix a principally polarized complex torus $(W,H,\Lambda)$
of complex dimension $g$. We now specialize the setting
of \S\ref{subsubsection:Setting} as follows:
$$V=W/\Lambda,\;\;\;\tilde{V}=W,\;\;\;G=\{w\mapsto \lambda+w\mid
\lambda\in\Lambda\}.$$
Also we put
$$\Phi(w):=H(w,w)$$
for all $w\in W$, thereby defining a Weil gauge. In the next
several paragraphs we recall  how  to classify and to construct
explicitly all
$\Phi$-normalized effective theta functions. We also recall the Riemann
quartic theta identity in a convenient form.

\subsubsection{Semicharacters} 
\label{subsubsection:Semicharacters}
A function $\psi$ on $\Lambda$ taking values in the group  of complex
numbers of absolute value $1$ is called a {\em semicharacter} of $\Lambda$
with
respect to
$H$ if 
$$\psi(\lambda+\mu)=\psi(\lambda)\psi(\mu)
\exp\left(\pi i \Im H(\lambda,\mu)\right)$$
for all $\lambda,\mu\in \Lambda$. The square of a semicharacter is a
unitary character, whence the terminology; also the ratio
of any two semicharacters is a unitary character. Real semicharacters
play an especially important role in the sequel, and these have the
following explicit description. Fix a symplectic
$\ZZ$-basis $\{A_i,B_i\}_{i=1}^g$ for $\Lambda$ arbitrarily. 
Let $A$ (resp., $B$) be the vector of length $g$ with entries $A_i$
(resp., $B_i$). Every real semicharacter
$\psi$ of $\Lambda$ with respect to $H$ takes the form
$$\psi(m\cdot A+n\cdot B)=(-1)^{m\cdot n+a\cdot m+b\cdot n}\;\;\;(m,n\in
\ZZ^g)$$ for some $a,b\in \ZZ^g$ uniquely determined modulo $2$.
It can be shown (for an indication of proof see the end of
\S\ref{subsubsection:ExplicitThetaConstruction} below) that the parity of
the inner product
$a\cdot b$ depends only on the real semicharacter
$\psi$, not on the choice of symplectic $\ZZ$-basis
$\{A_i,B_i\}$. We define the  {\em parity} of
$\psi$ to be that of $a\cdot b$.

\subsubsection{Classification of $\Phi$-normalized
effective theta functions}
\label{subsubsection:ThetaMultipliers}
According to \cite[Chap.\ 6]{Weil}, for each
semicharacter
$\psi$ of $\Lambda$ with respect to $H$ there exists
a not-identically-vanishing entire function
$\vartheta$ on $W$ unique up to a nonzero constant
factor such that
$$\vartheta(w+\lambda)=
\psi(\lambda)\exp\left(\pi H(\lambda,w)+
\frac{\pi}{2}
H(\lambda,\lambda)\right)\vartheta(w)
$$
for all $w\in W$ and $\lambda\in
\Lambda$. (We briefly sketch in
\S\ref{subsubsection:ExplicitThetaConstruction} below the
calculation  by which existence and uniqueness 
are proved.) We call any such function
$\vartheta$ a {\em theta function} of {\em type}
$(W,H,\Lambda, \psi)$.
 Clearly a theta function  of type $(W,H,\Lambda,\psi)$ is a
$\Phi$-normalized effective theta function in the sense of
Theorem~\ref{Theorem:ThetaTool}. Conversely, every
$\Phi$-normalized effective theta function $\vartheta$ in the sense
of Theorem~\ref{Theorem:ThetaTool}
is a theta function of type
$(W,H,\Lambda,\psi)$ for a uniquely determined semicharacter
$\psi$.  

\subsubsection{Natural operations on theta functions}
\label{subsubsection:NaturalOperations}
Let a theta function
$\vartheta(w)$ of type
$(W,H,\Lambda,\psi)$ be given. Then $\vartheta(-w)$ is a theta function
of type $\left(W,H,\Lambda,\overline{\psi}\right)$. In particular,
if $\overline{\psi}=\psi$, then $\vartheta(\pm w)=\pm \vartheta(w)$.
In other words, if $\psi$ is real, then $\vartheta$ has a well defined
parity. Given another semicharacter
$\psi'$ of
$\Lambda$ with respect to $H$ and $t\in W$ such that
$$\psi'(\lambda)=\psi(\lambda) \exp(2\pi i\Im
H(\lambda,t))$$
for all $\lambda\in \Lambda$,
 then $$\exp(-\pi
H(t,w))\vartheta(w+t)$$ is a theta function of type
$(W,H,\Lambda,\psi')$. In other words, roughly speaking,
any given $\Phi$-normalized theta function gives rise to all others
by translation and adjustment by elementary nowhere-vanishing factors.

\subsubsection{Explicit
construction of theta functions}
\label{subsubsection:ExplicitThetaConstruction}
Fix a
symplectic $\ZZ$-basis
$\{A_i,B_i\}_{i=1}^g$ for $\Lambda$ and a semicharacter
$\psi$ of $\Lambda$ with respect to $H$.   
For simplicity we assume that:
\begin{itemize}
\item $W=\CC^g$, the latter viewed for
computational purposes as the space of column vectors
of length $g$ with complex entries.
\item
$A_i$ is the
$i^{th}$ column of the $g$ by $g$ identity matrix.
\end{itemize}
Let $\tau$ be the $g$ by $g$ matrix defined by the
following condition:
\begin{itemize}
\item $B_i$ is the
$i^{th}$ column of $\tau$.
\end{itemize}
In this situation necessarily:
\begin{itemize}
\item $\tau$ is symmetric with positive definite imaginary part.
\item $\bar{v}^T(\Im \tau)^{-1}w=H(v,w)$
for all $v, w\in \CC^g$.
\end{itemize}
Moreover, there exist column
vectors
$a,b\in
\RR^g$ unique modulo $\ZZ^g$ with the following property:
\begin{itemize} 
\item $\psi(m +\tau n)=\exp\left(\pi i\,m^Tn+2\pi
i(m^T a-n^Tb)\right)$ for all $m,n\in
\ZZ^g$. 
\end{itemize}
As usual, cf.\ \cite[p.\ 1]{Fay}
or \cite[p.\ 123]{MumfordTataI},  put
$$\theta\left[\begin{array}{c}a\\b\end{array}\right]
(w,\tau ):=
\sum_{\ell\in \ZZ^g}
\exp\left(\pi i(\ell+a)^T\tau (\ell+a)+2\pi
i(\ell+a)^T(w+b)\right)
$$
thereby defining 
a holomorphic function of $w\in \CC^g$
that does not vanish identically.
To abbreviate we now drop reference to $a$, $b$ and $\tau$
since these are being held fixed, and we simply write
$\theta(w)$.  We have
$$\theta(w+m+\tau n)=
\exp\left(2\pi i (a^Tm-b^Tn)-\pi in^T\tau n -2\pi in^Tw\right)
\theta(w)
$$
for all $m,n\in \ZZ^g$ and $w\in \CC^g$. It is easy to show by
the method of ``undetermined Fourier coefficients'' that the system
of functional equations above characterizes $\theta(w)$ uniquely
up to a nonzero constant factor. By
a straightforward calculation it can be verified
that the function
$$\exp\left(\frac{\pi}{2}
w^T (\Im \tau )^{-1} w\right)
\theta(w)
$$
is a theta function of type $(W,H,\Lambda,\psi)$ and moreover the
only such up to a nonzero constant factor. From the explicit
presentation of theta functions of type $(W,H,\Lambda,\psi)$ just
recalled it follows in particular that if $\psi$ is real, then
the parity of
$\psi$ as defined at the end of \S\ref{subsubsection:Semicharacters}
coincides with the parity of any theta function of type
$(W,H,\Lambda,\psi)$.

\begin{Proposition}\label{Proposition:HeisenbergApplication}
Fix a semicharacter $\psi$ on $\Lambda$ with respect to $H$
and a theta function $\vartheta$ of type $(W,H,\Lambda,\psi)$.
Let $\Lambda'\subset W$ be a cocompact discrete
subgroup of $W$ such that
$(W,H,\Lambda')$ is also a principally polarized abelian variety.
Assume further that
$\card\left(\frac{\Lambda+\Lambda'}{\Lambda\cap\Lambda'}\right)<\infty$.
Fix a semicharacter $\psi'$ of $\Lambda'$ with respect to $H$
agreeing with $\psi$ on $\Lambda\cap \Lambda'$
and a theta function $\vartheta'$ of type $(W,H,\Lambda',\psi')$.
 Fix a (necessarily finite) set
of representatives $L\subset \Lambda$ for the quotient
$\Lambda/(\Lambda\cap\Lambda')$.
Then there exists a unique family 
$\{C_\lambda\}_{\lambda\in L}$ of
complex constants such that 
$$\vartheta(w)=\sum_{\lambda\in L}
C_\lambda\exp(-\pi
H(\lambda,w))\vartheta'(\lambda+w)$$
for all $w\in W$.
Morever, none of the constants $C_\lambda$ vanish.
\end{Proposition}
\proof 
We bring in a powerful idea developed at length in \cite{MumfordTataIII}.
Put
$$\HH:=\{[\lambda,s]\mid \lambda\in W,\;\;s\in
\CC,\;\;|s|=1\},$$ and equip $\HH$ with a group law by the
rule
$$[\lambda,s][\mu,t]:=[\lambda+\mu,st\exp(\pi i\Im
H(\lambda,\mu))],$$ 
thereby constructing the {\em Heisenberg group} naturally associated to 
the pair $(W,H)$. The group
$\HH$ acts naturally on the space of entire functions defined on $W$
by the rule
$$([\lambda,s]f)(w):=s\exp\left(-\pi
H(\lambda,w)-\frac{\pi}{2}H(\lambda,\lambda)\right)f(w+\lambda).$$
By definition of a semicharacter, the map
$$\left(\lambda\mapsto
\left[\lambda,\psi(\lambda)\right]\right):\Lambda\rightarrow\HH
$$
is an injective group homomorphism. We denote the image of this map
by
$\HH(\Lambda,\psi)$. For any unitary character $\chi$ of
$\Lambda$, a theta
function of type
$(W,H,\Lambda,\psi\chi)$ is the same thing as a not-identically-vanishing
holomorphic function $\varphi$ on $W$ tranforming under the action of
$\HH(\Lambda,\psi)$ by the rule
$$[\lambda,\psi(\lambda)]\varphi=\chi(\lambda)\cdot\varphi.$$ Let
$\Theta$ be the space of holomorphic functions on $W$ fixed under the
action of the group
$$\HH(\Lambda,\psi)\cap\HH(\Lambda',\psi')=
\{[\lambda,\psi(\lambda)]\mid \lambda\in \Lambda\cap \Lambda'\}\subset
\HH.$$ Then, so we claim, $\Theta$ is a regular complex representation
of the finite abelian group 
$$\HH(\Lambda',\psi')/(\HH(\Lambda,\psi)\cap
\HH(\Lambda',\psi'))$$ the isotypical
components of which are permuted simply transitively by the finite abelian
group
$$\HH(\Lambda,\psi)/(\HH(\Lambda,\psi)\cap\HH(\Lambda',\psi')).$$
Since the proof of the claim is merely a recapitulation of themes from the
proof of the Stone-von Neumann theorem (see \cite[Thm.\ 1.2,
p.3]{MumfordTataIII}), we omit the details. The claim granted, the result
immediately follows.
  \qed

\begin{Corollary}\label{Corollary:RiemannQuarticIdentity}
Fix a real semicharacter $\psi_0$ of $\Lambda$ with respect to $H$. Fix a
theta function
$\vartheta_0$ of type
$(W,H,\Lambda,\psi_0)$. 
Fix a set of representatives
$0\in M\subset
\frac{1}{2}\Lambda$ for the quotient $\frac{1}{2}\Lambda/\Lambda$.
Put
$$\vartheta_\mu(w):=\exp(-\pi H(\mu,w))\vartheta_0(w+\mu)$$
for all $\mu\in M$ and $w\in W$.  
Put
$$T:=\frac{1}{2}\left[\begin{array}{rrrr}
1&1&1&1\\
1&1&-1&-1\\
1&-1&1&-1\\
1&-1&-1&1
\end{array}\right].$$
Then there exists
a unique family 
$\{C_\mu\}_{\mu\in M}$ of complex constants
such that
$$
\prod_{i=1}^4\vartheta_0\left(w_i\right)=
\sum_{\mu\in M}
C_\mu
\prod_{i=1}^4\vartheta_\mu\left(
\sum_{j=1}^4 T_{ij}w_j\right)
$$
for all $w_1,w_2,w_3,w_4\in W$. Moreover, none of the constants
$C_\mu$ vanish. (This result is a translation into Weil-style
language of Riemann's
quartic theta identity. For a classical-style presentation of the
latter, see \cite[p.\ 212]{MumfordTataI}.)
\end{Corollary}

\proof It is a tedious job but not an especially
difficult one to verify that
Proposition~\ref{Proposition:HeisenbergApplication} applies with
$$\left[\begin{array}{c}
W\\
W\\
W\\
W\end{array}\right],\;\;\;
\left[\begin{array}{cccc}
H\\
&H\\
&&H\\
&&&H\end{array}\right],\;\;\;T\left[\begin{array}{c}
\Lambda\\
\Lambda\\
\Lambda\\
\Lambda
\end{array}\right],\;\;\;
\left[\begin{array}{c}
\Lambda\\
\Lambda\\
\Lambda\\
\Lambda
\end{array}\right]$$
in place of $W$, $H$, $\Lambda$ and $\Lambda'$, respectively,
the semicharacters
$$\left[\begin{array}{c}
\lambda_1\\
\lambda_2\\
\lambda_3\\
\lambda_4
\end{array}\right]\mapsto
\prod_{i=1}^4\psi_0\left(\sum_{j=1}^4T_{ij}\lambda_j\right),\;\;\;\;
\left[\begin{array}{c}
\lambda_1\\
\lambda_2\\
\lambda_3\\
\lambda_4
\end{array}\right]\mapsto
\prod_{i=1}^4\psi_0(\lambda_i)
$$ 
in
place of
$\psi$ and
$\psi'$, respectively, the theta functions
$$
\left[\begin{array}{c}
w_1\\
w_2\\
w_3\\
w_4
\end{array}\right]\mapsto \prod_{i=1}^4\vartheta_0\left(\sum_{j=1}^4
T_{ij}w_j\right),\;\;\;
\left[\begin{array}{c}
w_1\\
w_2\\
w_3\\
w_4
\end{array}\right]\mapsto \prod_{i=1}^4\vartheta_0\left(w_i\right)
$$
in place of $\vartheta$ and $\vartheta'$, respectively, 
and the set
$$\left\{\left.\left[\begin{array}{c}\mu\\
\mu\\
\mu\\
\mu\end{array}\right]\right|\mu\in M\right\}$$
in place of $L$. Accordingly, there exists a  unique family
$\{C_\mu\}_{\mu\in M}$ of complex constants such that
$$\prod_{i=1}^4\vartheta_0\left(\sum_{j=1}^4T_{ij}w_j\right)=
\sum_{\mu\in M}
C_\mu
\prod_{i=1}^4\vartheta_\mu\left(w_i\right)
$$
for all $w_1,w_2,w_3,w_4\in W$ and moreover none of the constants
$C_\mu$ vanish. Since $T^2=1$,
this last identity is equivalent to the desired one. \qed

\subsection{Example: Compact Riemann
surfaces}
\subsubsection{Specialization of the setting}
Let $X$ be a compact Riemann surface of
genus $g>0$, let $\tilde{X}\rightarrow X$ be a
universal covering map, and put $\Gamma:=\Aut(\tilde{X}/X)$.
We now specialize the setting of
\S\ref{subsubsection:Setting} to  the case
$$V=X,\;\;\;\tilde{V}=\tilde{X},\;\;\;G=\Gamma.$$
It is convenient to fix a basepoint $\infty\in X$
and a lifting $\tilde{\infty}\in \tilde{X}$ thereof.
In the next several paragraphs we work out a Weil-style analytic
description of the Jacobian of
$X$  and then we 
study the multiplier systems associated to theta functions on
$\tilde{X}$ relative to $\Gamma$ representing divisors on $X$ of degree
$0$ and of degree $g$.

\subsubsection{Basic topological notation}
Given points $\tilde{P},\tilde{Q}\in \tilde{X}$,
we denote by $[\tilde{P}\rightarrow\tilde{Q}]$ a choice of path in
$X$ admitting a lifting to a path issuing from $\tilde{P}$ and
terminating at $\tilde{Q}$. 
Given $1$-cycles $c_1$ and $c_2$ on $X$ in general
position, let $\card(c_1\cap c_2)$ denote the
signed number of intersections of $c_1$ with $c_2$,
where, as usual, we count $+1$ where $c_2$ crosses $c_1$ from
right to left and $-1$
at the other crossings.   Whenever we speak of paths, loops,
cycles, chains, etc., it is understood that all such are sufficiently
differentiable to integrate over.

\begin{Proposition}\label{Proposition:Polarization} Let
$W$ be the
$\CC$-linear dual of the space of differentials of the first kind on $X$.
Let $\Lambda$ be the subgroup of $W$
consisting of $\CC$-linear functionals of the form 
$\left(\omega\mapsto \int_c\omega\right)$ for some $1$-cycle $c$ on $X$.
There exists a unique hermitian form 
$$H:W\times
W\rightarrow \CC$$ antilinear on the left such that 
for all $1$-cycles $c_1$ and $c_2$ on $X$ in
general position we have
$$\Im H\left(\left(\omega\mapsto \int_{c_1}\omega\right),
\left(\omega\mapsto
\int_{c_2}\omega\right)\right)=\card(c_1\cap c_2).
$$ 
The triple $(W,H,\Lambda)$ is a principally polarized
abelian variety. (The triple $(W,H,\Lambda)$ is the Jacobian of $X$
described in Weil-style language. In the sequel we work exclusively with
this version of the Jacobian.)
\end{Proposition}
\proof By Hodge theory the subgroup $\Lambda$ is cocompact and discrete.
Uniqueness of $H$ is clear. The intersection pairing on the
$1$-dimensional homology of $X$ is well known to be alternating and to
have unit pfaffian. We have only to prove that $H$ with the desired
property exists and is positive definite. We construct a candidate $H_0$
for $H$ as follows. To each $w\in W$ we associate a holomorphic $1$-form
$\zeta_w$ by the rule $w=\left(\omega\mapsto\frac{1}{2}\int
\omega\wedge
\bar{\zeta}_w\right)$. The map $w\mapsto \zeta_w$ identifies
$W$ in $\CC$-antilinear fashion with the space of holomorphic
$1$-forms on $X$. Put
$$H_0(v,w)=\frac{i}{2}\int \zeta_v\wedge \bar{\zeta}_w$$
for all $v,w\in W$, thereby defining a positive definite hermitian
form $H_0$ on
$W$ antilinear on the left.
 For
each
$1$-cycle $c$ on $X$
there exists by Poincar\'{e} duality and Hodge theory a
unique holomorphic $1$-form $\zeta_c$ on $X$
such that  
$$\int_c\alpha=\int\alpha\wedge
\Re\zeta_c=\frac{1}{2}\int\alpha\wedge \bar{\zeta}_c$$ for all
smooth closed
$1$-forms $\alpha$ on $X$. It follows that for all $1$-cycles $c_1$
and
$c_2$ on $X$ in general position we have
$$\begin{array}{cl}
&\displaystyle\Im H_0\left(\left(\omega\mapsto
\int_{c_1}\omega\right),
\left(\omega\mapsto \int_{c_2}\omega\right)\right)\\\\
=&\displaystyle\Im \left(
\frac{i}{2}\int
\zeta_{c_1}\wedge
\bar{\zeta}_{c_2}\right)=\int\Re\zeta_{c_1}\wedge
\Re\zeta_{c_2}=\card(c_1\cap c_2).
\end{array}$$ 
Thus our candidate $H_0$ has all the desired properties, i.~e.,
$H=H_0$.
\qed

\subsubsection{Convenient abuses of notation} 
Hereafter we often treat $1$-chains as though they were
points of $W$. More precisely, given a $1$-chain $c$, we often just
write
$c$ where more properly we should write, say,
$\left(\omega\mapsto\int_c\omega\right)$. Further in this
line we also often treat elements of $\Gamma$ as if they were points of
$W$. More precisely, given $\sigma\in \Gamma$, we often just write
$\sigma$ where more properly we should write, say,
$\left(\omega\mapsto
\int_{[\tilde{\infty}\rightarrow\sigma\tilde{\infty}]}\omega\right)$.
This saves a lot of writing and should not cause any confusion.
\begin{Proposition}\label{Proposition:VanishingFact}
Let a differential $\omega$ of the first kind on $X$ and $w\in W$
be given. If $\int_c\omega=H(c,w)$ for all $1$-cycles on $X$,
then $\omega=0$ and $w=0$. 
\end{Proposition}
\proof We continue in the setting of the proof of
Proposition~\ref{Proposition:Polarization}. We have
$$\int_c\bar{\omega}=H(w,c)=\frac{i}{2}\int \zeta_{w}\wedge \bar{\zeta}_c
=i\int \zeta_{w}\wedge\Re \zeta_c=
\int_c i\zeta_{w},$$
hence the $1$-forms $\bar{\omega}$ and $i\zeta_w$ have
the same periods and hence they are equal.
But the latter $1$-form is holomorphic and the former antiholomorphic.
Therefore both must vanish.
\qed

\begin{Proposition}\label{Proposition:UnitaryCharacters}
Let $D$ be a divisor of $X$ of degree zero.\\
\begin{enumerate}
\item There exists  a theta function $\vartheta$ on $\tilde{X}$
relative to $\Gamma$ unique up to a nonzero constant factor such that
$\vartheta$ represents
$D$ and tranforms  under the action of
$\Gamma$ by a unitary character. \\
\item With $\vartheta$ as above, we have
$$\sigma^*\vartheta=
\exp\left( 2\pi i \Im
H(\partial^{-1}D,\sigma)\right)\cdot\vartheta$$
for all $\sigma\in \Gamma$ and $1$-chains $\partial^{-1}D$ on $X$ with
boundary 
$D$. \\
\item Every unitary character of $\Gamma$ thus appears in
association with some divisor of $X$ of degree zero.\\
\end{enumerate}
(The proposition is a restatement of the theorems of Abel
and Jacobi in Weil-style language.)
\end{Proposition}
\proof Statement 1 is equivalent to the
classical fact that there exists a unique differential
$\xi$ of the third kind on $X$ with residual divisor $D$ and pure
imaginary periods. (Of course statement 1 is also a very special
case of Theorem~\ref{Theorem:ThetaTool}.) Statements 1 and 2
granted, statement 3 is proved by a well known argument we need
not repeat.  It remains only to prove statement 2, and this is
just a matter of translating from classical-style language to
Weil-style language. We take care with the details in order to
check signs and factors of $2$ and $\pi$. 

Clearly
$d\log
\vartheta$ is the lifting of
$\xi$ and we have
$$\sigma^*\vartheta=
\exp\left(\int_{[\tilde{\infty}\rightarrow\sigma\tilde{\infty}]}
\xi\right)\vartheta$$
for all $\sigma \in \Gamma$, where the loop
$[\tilde{\infty}\rightarrow\sigma\tilde{\infty}]$ is chosen to
avoid the support of $D$. Now arbitrarily fix a
$1$-chain
$c_D$ with boundary $D$.
In order
to prove statement (ii), it suffices to verify that
$$\int_c \xi\equiv 2\pi i \Im H(c_D,c)\bmod{2\pi i\ZZ}$$
for all loops $c$ on $X$ avoiding the support of $D$.

In the usual way cut $X$ open to form a $4g$-sided polygon
and construct a homology basis $\{A_i,B_i\}_{i=1}^g$
in the standard configuration---meaning in particular that
$$\card(A_i\cap A_j)=0=\card(B_i\cap B_j)=0,\;\;\;
\card(A_i\cap B_j)=\delta_{ij}$$
for $i,j=1,\dots,g$. We may assume without loss of generality that
$c_D$ is contained in the interior of the polygon, and that $c$ is one of
the
$A$'s and
$B$'s. According to the classical reciprocity law for differentials of
first and third kinds  (see \cite[p.\ 230]{GriffithsHarris}) we have
$$\sum_{i=1}^g\left(\left(\int_{A_i}\omega
\right)\left(\int_{B_i}\xi\right)
-\left(\int_{B_i}\omega\right)\left( \int_{A_i}\xi\right)\right)=
2\pi i\int_{c_D}\omega$$
for all differentials $\omega$ of the first kind on $X$. In other words,
we have an identity
$$\sum_{i=1}^g \left(\left(\frac{1}{2\pi i}\int_{B_i}\xi\right) A_i-
\left(\frac{1}{2\pi i}\int_{A_i}\xi\right) B_i\right)=c_D$$
holding in $W$. This last equation is enough to finish the proof.
\qed

\subsubsection{The Abel map and associated Weil gauge}
\label{subsubsection:AbelMap}
Put
$$z(\tilde{P}):=
\left(\omega\mapsto\int_{[\tilde{\infty}\rightarrow\tilde{P}]}
\omega\right)\in
W$$ for all $\tilde{P}\in\tilde{X}$, thereby
defining the {\em
Abel map}
$$z:\tilde{X}\rightarrow W$$
based at $\tilde{\infty}$. In calculations below we repeatedly exploit
the relation
$$\sigma^*z=z+\sigma$$
for all $\sigma\in \Gamma$; of course, on the right, by abuse of
notation,
$\sigma$ stands in for the linear functional
$$\left(\omega\mapsto\int_{[\tilde{\infty}\rightarrow
\sigma\tilde{\infty}]}\omega\right)\in W.$$  Now let 
$$\Phi(w):=H(w,w)$$ be the Weil gauge
with which the Jacobian $(W,H,\Lambda)$ is canonically equipped.
Clearly the pull-back $z^*\Phi$ is a Weil
gauge on
$\tilde{X}$. Another more direct description of $z^*\Phi$ can be given as
follows. Let
$\zeta_1,\dots,\zeta_g$ be any
$\CC$-basis for the space of differentials of the first kind
on $X$ orthonormalized by the condition
$$\frac{i}{2}\int\zeta_i\wedge \bar{\zeta}_j=\delta_{ij}.$$
Then we have 
$$(z^*\Phi)(\tilde{P})=\sum_{i=1}^g
\left|\int_{[\tilde{\infty}\rightarrow\tilde{P}]}\zeta_i\right|^2$$
for all $\tilde{P}\in \tilde{X}$. From the latter description of $\Phi$
it  follows directly that
$\frac{i}{2}\partial\barpartial z^*\Phi$ is the lifting
of the closed real positive $(1,1)$-form
$$\frac{i}{2}\sum_{i=1}^g \zeta_i\wedge \bar{\zeta}_i.$$
In turn it follows by
Proposition~\ref{Proposition:PoincareRecognition} that the function
$z^*\Phi$ is a gauge for any divisor of degree $g$.

\subsubsection{Association of a semicharacter to each divisor of
degree
$g-1$} Fix a divisor
$D$ of $X$ of degree
$g-1$. By Theorem~\ref{Theorem:ThetaTool} the divisor
$D+\infty$ is represented by a 
$z^*\Phi$-normalized theta function on $\tilde{X}$ relative to $\Gamma$,
say
$\vartheta$, unique up to a nonzero constant factor. It is easy to see
that the multiplier system determined by $\vartheta$ has to be the
pull-back via the Abel map $z$  of the multiplier system determined by
some
$\Phi$-normalized theta function on $W$ relative to $\Lambda$. The
upshot is that there exists a unique semicharacter $\psi_D$ of $\Lambda$
with respect to $H$ such that
$$\sigma^*\vartheta
=\psi_D(\sigma)\cdot
\exp\left(\pi H(\sigma,z)
+\frac{\pi}{2}H(\sigma,\sigma)
\right)\cdot\vartheta$$
for all $\sigma\in \Gamma$. 
\begin{Proposition} Notation as in the
paragraph above, the following hold:\\
\begin{enumerate}
\item The semicharacter $\psi_D$ is independent of the choice of the
basepoint $\infty$ and lifting $\tilde{\infty}$ thereof.\\
\item The
construction
$D\mapsto
\psi_D$ puts the classes of divisors of $X$ of degree $g-1$ in bijective
correspondence with the semicharacters of $\Lambda$ with respect to
$H$.\\
\end{enumerate}
(The proposition is a  Weil-style
description of the correspondence between theta
functions with characteristics and divisors of degree $g-1$.
See \cite[Chap.\ II, \S3]{MumfordTataI} for a classical-style
treatment of this correspondence.)
\end{Proposition}
\proof Statement 1 granted, statement 2 follows immediately
from \linebreak Proposition~\ref{Proposition:UnitaryCharacters}. We turn
now to the proof of statement 1. Fix
$P\in X$ and a lifting
$\tilde{P}\in
\tilde{X}$ arbitrarily. By a repetition of the arguments made above
there exists a theta function
$\vartheta_1$ unique up to a nonzero constant multiple and
a unique semicharacter
$\psi_1$ of
$\Lambda$ with respect to
$H$ such that $\vartheta_1$ represents the divisor $D+P$ and
transforms according to the rule
$$\sigma^*\vartheta_1
=\psi_1(\sigma)\cdot
\exp\left(\pi H(\sigma,z-z(\tilde{P}))
+\frac{\pi}{2}H(\sigma,\sigma)
\right)\cdot\vartheta_1$$
for all $\sigma\in \Gamma$.
Then the theta function
$$\vartheta_2=\vartheta_1/\vartheta\cdot
\exp\left(\pi
H\left([\tilde{\infty}\rightarrow\tilde{P}],
z-z(\tilde{P})\right)\right)$$
 represents the divisor
$P-\infty$ of degree zero and transforms according to the unitary rule
$$\sigma^*\vartheta_2=\left(\psi_1/\psi_D\right)(\sigma)
\cdot \exp\left(2\pi i\Im H\left([\tilde{\infty}\rightarrow
\tilde{P}],\sigma\right)\right)\cdot\vartheta_2$$
for all $\sigma\in \Gamma$.
 By Proposition~\ref{Proposition:UnitaryCharacters} it follows that
$\psi_1=\psi_D$.
\qed

\subsection{Example: The prime form}
\subsubsection{Specialization of the setting}
Fix a compact Riemann surface $X$ of genus
$g>0$ and a universal covering map $\tilde{X}\rightarrow X$.
Fix a
basepoint
$\infty\in X$ and lifting $\tilde{\infty}\in \tilde{X}$ thereof.
Put $\Gamma:=\Aut(\tilde{X}/X)$. 
 We specialize the setting of
\S\ref{subsubsection:Setting} to  the case
$$V=X\times X,\;\;\;\tilde{V}=\tilde{X}\times \tilde{X},\;\;\;
G=\Gamma\times \Gamma.$$
 Let 
$$\Delta\subset X\times X$$
be the diagonally embedded copy of $X$ and consider the divisor
$$\Delta'=-X\times \infty-\infty\times
X+\Delta.$$ 
We are going to find a natural choice of gauge for $\Delta'$ and then
we are going to calculate the multiplier system of the correspondingly
normalized theta function representing
$\Delta'$.

\subsubsection{A gauge for $\Delta'$}
\label{subsubsection:EGauge}
Let
$(W,H,\Lambda)$ be the Jacobian of $X$. Let 
$$\Phi(w):=H(w,w)$$
be the Weil gauge naturally associated to the Jacobian.
Let $z$ denote
the Abel map $\tilde{X}\rightarrow W$ based at $\tilde{\infty}$.
For $i=1,2$, let $z^{(i)}$ denote the $i^{th}$ projection
$\tilde{X}\times \tilde{X}\rightarrow\tilde{X}$ followed by the Abel map
$z$. 
Put\\
$$
\Psi:=
\Phi\left(z^{(1)}-z^{(2)}\right)-
\Phi\left(z^{(1)}\right)-\Phi\left(-z^{(2)}\right).$$
thereby defining a Weil gauge $\Psi$ on $\tilde{X}\times \tilde{X}$.
 We claim that $\Psi$ is a gauge for $\Delta'$.
Let
$$p_1,p_2:X\times X\rightarrow X$$ be the two projections
and let
$\zeta_1,\dots,\zeta_g$ be a basis for the holomorphic
$1$-forms on $X$ orthonormalized by the condition
$$\frac{i}{2}\int\zeta_i\wedge \bar{\zeta}_j=\delta_{ij}.$$ 
Then $\frac{i}{2}\partial\barpartial\Psi$
is the lifting of the real closed
$(1,1)$-form
$$\begin{array}{rcl}
\alpha&=&\displaystyle\frac{i}{2}\sum_{i=1}^g
\left((p_1^*\zeta_i-p_2^*\zeta_i)\wedge
(p_1^*\bar{\zeta}_i-p_2^*\bar{\zeta}_i)-p_1^*(\zeta_i\wedge
\bar{\zeta}_i)-p_2^*((-\zeta_i)\wedge
(-\bar{\zeta}_i))\right)\\\\
&=&\displaystyle \frac{i}{2}\sum_{i=1}^g
\left(-p_1^*\zeta_i\wedge p_2^*\bar{\zeta}_i+
p_1^*\bar{\zeta}_i\wedge p_2^*\zeta_i\right)
\end{array}
$$ on $X$.  By
Proposition~\ref{Proposition:PoincareRecognition} the proof of the
claim boils down to  verifying the identity
$$\int_{X\times X}\alpha\wedge \beta=
\int_\Delta\beta-\int_{\infty\times X}\beta-\int_{X\times
\infty}\beta$$ for $\beta$ ranging over the $\CC$-basis
$$\frac{i}{2}p_1^*\zeta_i\wedge p_2^*\zeta_j,\;\;\;
\frac{i}{2}p_1^*\bar{\zeta}_i\wedge
p_2^*\bar{\zeta}_j,\;\;\;\frac{i}{2}p_1^*\zeta_i\wedge
p_2^*\bar{\zeta}_j,\;\;\; \frac{i}{2}p_1^*\bar{\zeta}_i\wedge
p_2^*\zeta_j$$
and
$$\frac{i}{2}p_1^*(\zeta_1\wedge \bar{\zeta}_1),\;\;\;
\frac{i}{2}p_2^*(\zeta_1\wedge \bar{\zeta}_1)$$
for the de Rham cohomology of $X\times X$ in dimension
$2$. The latter calculation is straightforward and
can safely be omitted.  Thus the claim is proved.

\subsubsection{The prime form $E$}
We define the
{\em prime form} $E$ to be  the \linebreak $\Psi$-normalized
theta function on $\tilde{X}\times
\tilde{X}$ unique up to a nonzero constant
factor representing the divisor $\Delta'$. The existence of $E$ is
guaranteed by Theorem~\ref{Theorem:ThetaTool}, and the uniqueness of
$E$ is clear. The notion of prime form defined here is nearly but not
exactly the same as the notion considered in
\cite[Chap.\ 2]{Fay} or \cite[pp.~3.207-3.213]{MumfordTataII}. We
omit discussion of the comparison, just remarking that 
\cite[Lemma 2, Chap.\ IIIb,
\S1, p.\ 3.211]{MumfordTataII} can be used to make the prime form
as defined here explicit. In any case an explicit formula will not be
needed. For
our purposes it suffices simply to know that
$E$ exists, is unique up to a nonzero constant factor, and has the
transformation properties summarized in 
Proposition~\ref{Proposition:PrimeTransformationProperties}
below.

\subsubsection{Guessing the multiplier system for $E$}
The form in which we presented the definition of the gauge $\Psi$ was
intended to suggest the following procedure for guessing the
multiplier system of $E$. Fix a semicharacter
$\psi$ of
$\Lambda$ with respect to $H$ arbitrarily. Consider the multiplier
system
$$\left\{\psi(\lambda)\exp\left(\pi
H(\lambda,w)+\frac{\pi}{2}H(\lambda,\lambda)\right)\right\}$$
on $W$ relative to $\Lambda$ determined by a theta function of type
$(W,H,\Lambda,\psi)$.  Pulling back under the map $z^{(1)}$ we
obtain a multiplier system
$$\left\{\psi(\sigma_1)\exp\left(\pi H\left(\sigma_1,z^{(1)}\right)
+\frac{\pi}{2}H(\sigma_1,\sigma_1)\right)\right\}$$
on $\tilde{X}\times \tilde{X}$ relative to $\Gamma\times \Gamma$.
Similarly, pulling back under the map $-z^{(2)}$ we obtain a multiplier
system
$$\left\{\psi(-\sigma_2)\exp\left(\pi H\left(-\sigma_2,-z^{(2)}\right)+
\frac{\pi}{2}H(-\sigma_2,-\sigma_2)\right)\right\},$$
and pulling back under the map $z^{(1)}-z^{(2)}$ we obtain a multiplier
system
$$
\left\{
\psi(\sigma_1-\sigma_2)\exp\left(\pi
H\left(\sigma_1-\sigma_2,z^{(1)}-z^{(2)}\right)
+\frac{\pi}{2}H\left(\sigma_1-\sigma_2,\sigma_1-\sigma_2\right)\right)
\right\}.
$$
After dividing the last multiplier system by the product of the first two
and simplifying, we obtain a multiplier system
$$\left\{\exp\left(-\pi \left(H\left(\sigma_1,z^{(2)}\right)
+H\left(\sigma_2,z^{(1)}\right)+H(\sigma_1,\sigma_2)\right)\right)
\right\}$$
on $\tilde{X}\times \tilde{X}$ relative to $\Gamma\times \Gamma$ such
that any theta function determining that multiplier system is necessarily
$\Psi$-normalized. This is our guess for the multiplier
system determined by
$E$. Next, we prove the guess.

\begin{Proposition}\label{Proposition:PrimeTransformationProperties}
We have
$$(\sigma_1,\sigma_2)^*E=\exp\left(-\pi \left(
H(\sigma_1,z^{(2)})+ H(\sigma_2, z^{(1)})+
H(\sigma_1,\sigma_2)\right)\right)\cdot E$$
for all 
$$(\sigma_1,\sigma_2)\in \Gamma\times \Gamma=\Aut(\tilde{X}^2/X^2),$$
cf.\
\cite[p.\ 3.210]{MumfordTataII}.
Moreover, $E$ is
antisymmetric under exchange of factors in the product
$\tilde{X}\times \tilde{X}$.
\end{Proposition}
\proof The proof is in essence an adaptation to the present situation of
the proof of the theorem of the square. We introduce the following temporary
notation. We write
$(\sigma,\tau)$ instead of $(\sigma_1,\sigma_2)$,
we denote the multiplier system claimed for $E$ 
by 
$\{F_{\sigma\tau}\}$, and we denote the actual multiplier system of $E$
by  $\{F'_{\sigma\tau}\}$. 
The ratio
$\{F_{\sigma\tau}/F'_{\sigma\tau}\}$ is a unitary and hence
constant system of multipliers.  The restricted multiplier
systems
$$\left\{F_{\sigma 1}\mid_{\tilde{X}\times
\tilde{\infty}}\right\}_{\sigma\in
\Aut(\tilde{X}/X)},
\;\;\;\left\{F_{1\tau}\mid_{
\tilde{\infty}\times \tilde{X}}\right\}_{\tau\in
\Aut(\tilde{X}/X)}
$$
are identically equal to $1$, and the analogous remark holds
for
$\{F'_{\sigma\tau}\}$ since the invertible sheaf ${\mathcal O}_{X\times
X}(\Delta')$ has trivial restrictions to \linebreak $X\times \infty$ and
$\infty\times X$. Therefore the multiplier system
$\{F_{\sigma\tau}/F'_{\sigma\tau}\}$ is identically equal to $1$,
i.~e., the prime form $E$ transforms in the claimed fashion under the
action of $\Gamma\times \Gamma$. 
Since exchange of factors in the product $\tilde{X}\times \tilde{X}$
preserves both
$\Psi$ and
$\Delta'$, and hence can alter $E$ only by a nonzero constant factor,
$E$ is either symmetric or antisymmetric. The
sign is nailed down by considering what happens near the diagonal. 
\qed

\begin{Corollary}\label{Corollary:Conjugation}
For all canonical divisors $K$ and divisors $D$ of degree $g-1$
on $X$ we have 
$$\overline{\psi}_D=\psi_{K-D}.$$
(In particular, the real semicharacters of
$\Lambda$ with respect to $H$ are in canonical bijective correspondence
with the half-canonical divisor classes on $X$, cf.\ \cite[Chap.\ IIIa,
\S6]{MumfordTataII}.)
\end{Corollary}
\proof Let 
$$\delta:\tilde{X}\rightarrow \tilde{X}\times \tilde{X}$$
be the diagonal mapping. Note that 
$$\delta^*\Psi=-2z^*\Phi.$$
Let
$\{F_{\sigma\tau}\}$ temporarily denote the multiplier system determined
by
$E$.  
On the one hand, 
since the diagonal restriction of the invertible sheaf
${\mathcal O}_{X\times X}(-\Delta')$ is isomorphic to
$\Omega_X(2\infty)$,  
 it follows in view of the remark of
\S\ref{subsubsection:BundleGauge} that the diagonally restricted
 multiplier system
$$\left\{\delta^*F_{\sigma\sigma}\right\}_{\sigma\in
\Gamma}=
\{\exp(-2\pi H(\sigma,z)-\pi
H(\sigma,\sigma))\}_{\sigma\in\Gamma}
$$
is that determined by a $-2z^*\Phi$-normalized theta function on
$\tilde{X}$ relative to $\Gamma$ representing
$-K-2\infty$. On the other hand, since the reciprocal of the diagonally
restricted multiplier system can be written as the product of
multiplier systems associated to
$z^*\Phi$-normalized theta functions on $\tilde{X}$ relative
to $\Gamma$ representing
$D+\infty$ and
$K-D+\infty$, respectively, it is clear that
$\psi_D\psi_{K-D}=1$.
\qed

\subsection{Example: Determinant identities satisfied by the Riemann
theta function}
\label{subsection:DeterminantIdentities}

\subsubsection{Setting and notation}\label{subsubsection:ThetaSetting}
We now combine elements of all the preceding examples.
Fix a compact Riemann
surface
$X$ of genus
$g>0$ and a universal covering map $\tilde{X}\rightarrow X$. Put
$\Gamma:=\Aut(\tilde{X}/X)$.
 Fix a basepoint
$\infty\in X$ and lifting
$\tilde{\infty}\in \tilde{X}$ thereof.
 Let
$(W,H,\Lambda)$ be the Jacobian of $X$.
 Let
$z$ be the Abel map $\tilde{X}\rightarrow W$ based at
$\tilde{\infty}$.
 Fix an integer $n\geq g$.
 Given any function $f$ defined on $\tilde{X}$ and index
 $i=1,\dots,n$, let $f^{(i)}$  denote the result of following
the
$i^{th}$ projection 
$\tilde{X}^n\rightarrow \tilde{X}$ by $f$.
Let $E$ be the prime form. 
Given  indices $1\leq i<j\leq n$, let $E^{(i,j)}$ denote the result
of following the $(i,j)^{th}$ projection 
 $\tilde{X}^n\rightarrow
\tilde{X}\times \tilde{X}$ by $E$. 
The symbol $\propto$ placed between
two expressions indicates that both expressions define
not-identically-vanishing meromorphic functions on the same complex
manifold agreeing up to a nonzero constant factor.
 
\subsubsection{A special class of theta functions}
\label{subsubsection:SpecialThetaClass}
Given a semicharacter $\psi$ on $\Lambda$ with
respect to $H$, an integer $\ell$, and a meromorphic function
$u$ on $\tilde{X}$, we say that $u$ is a {\em theta function} of
{\em type}
$(X,\infty,\psi,\ell)$ if 
$$\sigma^*u=\psi(\sigma) \cdot\exp\left(\pi
H(\sigma,z)+\frac{\pi}{2}H(\sigma,\sigma)\right)\cdot
u$$
for all $\sigma\in G$ and morever $u$ has no singularities save 
poles of order at most $\ell-1$
at each lifting of the basepoint $\infty$. 
The
collection of all such theta functions $u$ forms a vector space over
the complex numbers.  
(Deviating from the convention we have followed above, in this case we
do not exclude the case $u=0$.)
Note that for any divisor $D$ of degree
$g-1$, integer $\ell$ and
$z^*\Phi$-normalized theta function $\phi$ representing
$D+\infty$ the map
$$\left(f\mapsto (\mbox{lifting of $f$})\cdot\phi
\right)$$
$$:H^0(X,\OO_X(D+\ell\infty))
\rightarrow\left(\begin{array}{l}
\mbox{space of theta functions}\\
\mbox{of type $(X,\infty,\psi_D,\ell)$}
\end{array}\right)$$ 
is bijective. Since every semicharacter $\psi$ is of the form $\psi_D$
for some $D$, and since for any such $D$ we have by assumption
$$\deg (D+n\infty)\geq g-1+g>2g-2,$$ it follows that
$$\dim_\CC\left(\begin{array}{l}
\mbox{space of theta functions}\\
\mbox{of type
$(X,\infty,\psi,n)$}
\end{array}\right)=n$$
 by Riemann-Roch.

\begin{Lemma}\label{Lemma:Characterization}  
Fix a semicharacter $\psi$ on $\Lambda$ with respect to $H$.
Recall that we are assuming that $n\geq g$. There
exists a not-identically-vanishing meromorphic function $\phi$ on
$\tilde{X}^n$ unique up to a nonzero constant factor  with the following
properties:\\
\begin{enumerate}
\item The function $\phi$ transforms according to the
rule
$$\sigma^*\phi=
\prod_{i=1}^n \psi(\sigma_i)\exp\left(\pi H\left(\sigma_i,z^{(i)}\right)\right)
\cdot \phi$$
for all $$\sigma=\left(\sigma_i\right)_{i=1}^n\in 
\Gamma^n=\Aut(\tilde{X}^n/X^n).$$ In particular, $\phi$ is a
theta function on $\tilde{X}^n$ relative to $\Gamma^n$.\\
\item The divisor 
$$-(n-1)\left(\infty\times X^{n-1}+X\times \infty\times X^{n-2}
+\cdots+X^{n-1}\times \infty\right)$$ 
is a lower bound for the divisor of $X^n$ represented by 
$\phi$. In particular, $\phi$ is regular on
$(\tilde{X}\setminus\Gamma\tilde{\infty})^n$.\\
\item The function $\phi$ is anti-symmetric under exchange of
factors of the product
$\tilde{X}^n$.\\
\end{enumerate}
\end{Lemma}
\proof 
Let $T\subset \tilde{X}\setminus
\Gamma\tilde{\infty}$ be a set of cardinality $n$ such that
the map
$$\left(u\mapsto u\mid_T\right):
\left(\begin{array}{l}
\mbox{space of theta}\\
\mbox{functions of}\\
\mbox{type $(X,\infty,\psi,n)$}
\end{array}\right)\rightarrow
\left(\begin{array}{l}
\mbox{space of}\\
\mbox{functions on $T$}
\end{array}\right)
$$ map is bijective.
 For $i=1,\dots,n$, by induction on $i$,
a meromorphic function on $\tilde{X}^n$ satisfying
conditions 1 and 2 vanishes identically if
it vanishes 
identically on $T^i\times \tilde{X}^{n-i}$. Thus we find that the
map
$$\left(\phi\mapsto \phi\mid_{T^n}\right):
\left(\begin{array}{l}
\mbox{space of functions}\\
\mbox{satisfying 1 and 2}
\end{array}\right)\rightarrow
\left(\begin{array}{l}
\mbox{space of}\\
\mbox{functions on $T^n$}
\end{array}\right)
$$
is bijective. Under
the latter map the space of functions satisfying conditions 1, 2 and 3
maps bijectively to the space of functions on $T^n$ antisymmetric under
exchange of factors. Clearly the latter space is one-dimensional.
\qed

\begin{Proposition}
\label{Proposition:BasicDeterminantIdentity}
Let $\vartheta$ be a theta function of type
$(W,H,\Lambda,\psi)$. Recall that we are assuming that $n\geq g$. Let
$u_1,\dots,u_n$ be a
$\CC$-basis for the space of theta functions of type
$(X,\infty,\psi,n)$. We have
$$\det_{i,j=1}^n u^{(i)}_j\propto
\vartheta\left(\sum_{i=1}^n
z^{(i)}\right)\prod_{1\leq i<j\leq n}E^{(i,j)}.$$
(The
proposition is a Weil-style formulation of a classical determinant
identity satisfied by the Riemann theta function. See
\cite[Prop.\ 2.16, p.\ 29]{Fay} for a classical-style formulation of
this identity.)
\end{Proposition}
\proof The
expression on the left (resp., right) side of the claimed relation
 defines a not-identically-vanishing meromorphic function on
$\tilde{X}^n$, and moreover this function clearly satisfies conditions
1, 2 and 3 (resp., 2 and 3)  enunciated in
Lemma~\ref{Lemma:Characterization}. It remains only to verify that the
function on the right satisfies condition 1. Fix
$$\sigma=(\sigma_i)_{i=1}^n\in \Gamma^n=\Aut(\tilde{X}^n/X^n)$$
arbitrarily. We have
$$\begin{array}{cl}
&\sigma^*\vartheta\left(\sum z^{(i)}\right)\\
=&
\vartheta\left((\sum z^{(i)})+(\sum\sigma_i)\right)\\
=&\psi(\sum \sigma_i)\exp(\pi H(\sum\sigma_i,\sum z^{(i)})+
\frac{\pi}{2}H(\sum\sigma_i,\sum\sigma_i))\vartheta(\sum z^{(i)}),\\\\
&\sigma^* E^{(i,j)}\\
=&((\sigma_i,\sigma_j)^*E)^{(i,j)}\\
=&
\exp(-\pi(H(\sigma_i,z^{(j)})+H(\sigma_j,z^{(i)})+H(\sigma_i,\sigma_j)))
E^{(i,j)}
\end{array}$$
by definition of a theta function of type $(W,H,\Lambda,\psi)$
and the
transformation
law for the
prime form
$E$
enunciated
in Proposition~\ref{Proposition:PrimeTransformationProperties},
respectively.
This is enough to finish the
proof. 
\qed

\begin{Corollary}
\label{Corollary:WRVT}
For any divisor
$D$ of
$X$ of degree
$g-1$, theta function
$\vartheta$ of type
$(W,H,\Lambda,\psi_D)$, divisor $D'$ of $X$ of degree $0$
and $1$-chain $\partial^{-1}D'$ on $X$ with boundary $D'$ 
we have
$$\vartheta\left(\partial^{-1}D'\right)=0
\Leftrightarrow h^0\left(D-D'\right)>0.
$$
(The proposition is a Weil-style
formulation of a weak version of the Riemann Vanishing Theorem.
For classical-style discussions of the latter, see 
\cite[Chap.\
II,
\S3]{MumfordTataI} or \cite[Riemann-Kempf Singularity
Theorem, p.\ 348]{GriffithsHarris}.)
\end{Corollary}
\proof Both sides of the logical equivalence to be proved depend only
on the divisor class of $D'$. We may therefore assume without loss of
generality that for some
$n\geq g$, distinct points 
$P_1,\dots,P_n\in X\setminus\infty$
and corresponding liftings
$\tilde{P}_1,\dots,\tilde{P}_n\in \tilde{X}$ we
have
$$D'=-n\infty+\sum_{i=1}^n P_i,\;\;\;
\partial^{-1}D'=\sum_{i=1}^n [\tilde{\infty}\rightarrow
\tilde{P}_i]=\sum_{i=1}^n z(\tilde{P}_i)\in W.$$ 
We then have
$$h^0\left(D+n\infty-\sum_{i=1}^n
P_i\right)=\dim_\CC\left(\begin{array}{l}
\mbox{space of theta functions}\\
\mbox{of type $(X,\infty,\psi,n)$}\\
\mbox{vanishing at $\tilde{P}_1,\dots,\tilde{P}_n$}
\end{array}\right)$$
by the remarks of \S\ref{subsubsection:SpecialThetaClass} and hence
$$h^0\left(D+n\infty-\sum_{i=1}^n
P_i\right)>0\Leftrightarrow\vartheta\left(\sum_{i=1}^n
z(\tilde{P}_i)\right)=0$$ by
Proposition~\ref{Proposition:BasicDeterminantIdentity}.
\qed

\subsubsection{Remark}
We continue in the setting of Corollary~\ref{Corollary:WRVT}.
At full strength, the Riemann Vanishing Theorem says that
$h^0(D-D')$ equals the order of vanishing of $\vartheta$
at the point $\partial^{-1}D'\in W$. For an extended discussion of the
Riemann Vanishing Theorem from the classical point of view, see, e.~g.,
\cite[loc.\ cit.]{GriffithsHarris}.
One can deduce the Riemann Vanishing Theorem at full strength
from Proposition~\ref{Proposition:BasicDeterminantIdentity} by
manipulation of symmetric functions on $X^n$. Such manipulations are
standard in soliton theory, cf.\
\cite[pp.\ 48-52]{SegalWilson}.

\begin{Corollary}\label{Corollary:PrimeTrick}
Fix a half-canonical divisor $D$ on $X$ 
and a theta function $\vartheta$ of type $(W,H,\Lambda,\psi_D)$.
The following are
equivalent:
\begin{enumerate}
\item $\vartheta$ vanishes at the origin to order at
least
$2$.
\item $h^0(D)\geq 2$.
\item The function
$\vartheta\left(z^{(1)}-z^{(2)}\right)$ on
$\tilde{X}\times \tilde{X}$ vanishes identically.
\end{enumerate}
(To prove the equivalence of conditions 2 and 3 only the
``weak'' Riemann Vanishing Theorem is needed.)
\end{Corollary}
\proof (1 $\Leftrightarrow$ 2) This equivalence is a special case of
the ``full strength'' Riemann Vanishing Theorem. 

(2 $\Rightarrow$ 3) If $h^0(D)>1$, then $h^0(D-P)>0$ for
all points
$P$ of $X$ and hence $h^0(D-P+Q)>0$ for all points $P$ and $Q$ of
$X$, whence the result 
via 
Corollary~\ref{Corollary:WRVT}.

(not 2 $\Rightarrow$ not 3) We may assume without loss of generality
that
$h^1(D)=1$, for otherwise $\vartheta(0)\neq 0$ by
Corollary~\ref{Corollary:WRVT} and there is nothing to prove. To finish
the proof we now follow the proof of
\cite[Lemma 2, p.\ 3.211]{MumfordTataII}. Pick a point $P$
on $X$ arbitrarily. We have either
$h^0(D+P)=1$
or else, by Riemann-Roch,
$$h^0(D-P)=h^0(K-D-P)=h^0(D+P)-1=1.$$ 
Existence of $Q$ such that $h^0(D+P-Q)=0$ 
in the former case and $h^0(D-P+Q)=0$ in the latter case is clear,
whence in both cases
the result via Corollary~\ref{Corollary:WRVT}.
\qed

\subsection{Example: An {\em ad hoc} deformation theory}
We continue to work in the setting of
 \S\ref{subsection:DeterminantIdentities}.

\begin{Proposition}\label{Proposition:Frame}
Fix a divisor $D$ of $X$ of degree $g-1$.
Fix a theta function $\vartheta$ of type
$(W,H,\Lambda,\psi_D)$. Make now the stronger assumption that $n\geq g+1$.
Fix distinct points 
$$P_1,\dots,P_n\in
X\setminus\{\infty\}$$
and corresponding liftings
$$\tilde{P}_1,\dots,\tilde{P}_n\in
\tilde{X}$$ such that 
$$\vartheta\left(\sum_{i=1}^n
z(\tilde{P}_i)\right)\neq 0.$$ For $i=1,\dots,n$ and $\tilde{P}\in
\tilde{X}\setminus\Gamma\tilde{\infty}$, put
$$u_i(\tilde{P}):=\vartheta\left(z(\tilde{P})+\sum_{\alpha\in
\{1,\dots,n\}\setminus\{i\}}
z(\tilde{P}_\alpha)\right)\prod_{\alpha\in
\{1,\dots,n\}\setminus\{i\}}
E\left(\tilde{P},\tilde{P}_\alpha\right),$$ thereby defining a
function regular on $\tilde{X}\setminus\Gamma\tilde{\infty}$
and meromorphic on $\tilde{X}$.
 The
following hold:
\begin{enumerate}
\item The function $u_i$ is a theta
function of type $(X,\infty,\psi_D,n)$. 
\item The functions
$u_1,\dots,u_n$ form a $\CC$-basis
for the space of theta functions of type $(X,\infty,\psi_D,n)$.
\item The divisor of $X$ represented by $u_i$, say $D_i$,
belongs to the divisor class of
$D+\infty$.  
\item The divisors $D_1+(n-1)\infty,\dots,D_n+(n-1)\infty$ are effective
and have no zeroes in common.
\end{enumerate}
\end{Proposition}
\proof 1. This is clear in view of our complete knowledge
of the laws governing the transformation of $\vartheta$ and $E$.

2. By construction of the $u_i$ we have
$$\prod_{i=1}^n u_i(\tilde{P}_i)=\det_{i,j=1}^n
u_i(\tilde{P}_j)\neq 0$$
and hence the $u_i$ are $\CC$-linearly independent.

3. By definition of $\psi_D$, a
$z^*\Phi$-normalized theta function representing
$D+\infty$ transforms under the action of $\Gamma$ according to the same
law as does any theta function of type $(X,\infty,\psi_D,n)$.

4. It is clear that the divisors $D_i+(n-1)\infty$ are effective. But
further, under our assumption that
$n\geq g+1$ we have
$$\deg ( D+n\infty)\geq (g-1)+(g+1)=
2g,$$ 
and hence the linear system of effective divisors belonging to the
divisor class of
$D+n\infty$ is basepoint free. 
\qed

\begin{Proposition}\label{Proposition:DeformationFormula}
We continue in the setting of the proposition above.
Fix $w_0\in W$. Put
$$\dot{\vartheta}(w):=\left.
\frac{\partial}{\partial s}\log
\vartheta\left(w+sw_0\right)\right|_{s=0},$$
thereby defining a meromorphic function on $W$ regular away from the zero
locus of $\vartheta$. For $i=1,\dots,n$ put
$$\dot{u}_i:=\dot{\vartheta}\left(z+\sum_{\alpha\in
\{1,\dots,n\}\setminus\{i\}}z(\tilde{P}_\alpha)\right),$$
thereby defining a meromorphic function on $\tilde{X}$ regular away from the
lifting to $\tilde{X}$ of the support of the effective divisor
$D_i+(n-1)\infty$.  If the functions
$\dot{u}_i$ differ from one another by constants,  then 
each of the functions $\dot{u}_i$ reduces to a constant and $w_0=0$.
\end{Proposition}
\proof
For $i=1,\dots,n$, let $U_i$ be the complement of the support of the
effective divisor $D_i+(n-1)\infty$. 
By Proposition~\ref{Proposition:Frame} the family 
$\{U_i\}$ is an open covering of $X$.
Let
$\tilde{U}_i\subset \tilde{X}$ be the inverse image of $U_i$ under the
covering map. Since $\dot{u}_i$ is regular on $\tilde{U}_i$
and the family
$\{\tilde{U}_i\}$ covers
$\tilde{X}$, it follows that there exists some holomorphic function
$\dot{u}$ on $\tilde{X}$ from which
each function
$\dot{u}_i$ differs by a constant. 
Clearly we have
$$\dot{\vartheta}(w+\lambda)=
\pi H(\lambda,w_0)+\dot{\vartheta}(w)$$
for all $\lambda\in \Lambda$ and $w\in W$ such that $\vartheta(w)\neq 0$
and hence
$$\sigma^*\dot{u}=\pi H(\sigma,w_0)+\dot{u}$$
for all $\sigma\in \Gamma$.
It follows that $\dot{u}$ is the primitive
of some holomorphic $1$-form on $X$, say $\omega$. By
construction of
$\omega$ we have
$$\int_{[\tilde{\infty}\rightarrow\sigma\tilde{\infty}]}
\omega=\pi H(\sigma,w_0)$$
for all $\sigma\in \Gamma$
and hence $\omega=0$ and $w_0=0$ by
Proposition~\ref{Proposition:VanishingFact}.
\qed

\section{The edited $4\Theta$-embedding of
a Jacobian}\label{section:AbeliantIdentity}
\subsection{The edited $4\Theta$ linear system}
\subsubsection{The setting}
Fix a
compact Riemann surface
$X$ of genus
$g>0$ and a universal covering map $\tilde{X}\rightarrow X$.
Put $\Gamma:=\Aut(\tilde{X}/X)$.
 Fix a basepoint
$\infty\in X$ and lifting
$\tilde{\infty}\in \tilde{X}$ thereof.
 Let
$(W,H,\Lambda)$ be the Jacobian of $X$.
 Let
$z$ be the Abel map $\tilde{X}\rightarrow W$ based at
$\tilde{\infty}$.
 Let $E$ be the prime form. 
Fix a real semicharacter $\psi_0$ of $\Lambda$ with respect to $H$
and a theta function $\vartheta_0$ of type $(W,H,\Lambda,\psi_0)$.
More generally, put
$$\vartheta_t(w)=\exp\left(-\pi
H(t,w)\right)\vartheta_0(w+t)$$ for all $t\in W$ and
correspondingly put
$$\psi_t(\lambda)=\exp(2\pi i\Im
H(\lambda,t))\psi_0(\lambda)$$
for all $\lambda \in\Lambda$,
so that then $\vartheta_t$ is a theta function of type
$(W,H,\Lambda,\psi_t)$. Fix a set
$$0\in M\subset \frac{1}{2}\Lambda$$ of
representatives for the quotient
$\frac{1}{2}\Lambda/\Lambda$. Then
$$\{\psi_\mu\}_{\mu\in M}$$
is
the family of $4^g$ real semicharacters of $\Lambda$ with respect to
$H$.

\subsubsection{Definition}
By the Lefschetz embedding theorem the family of $4^g$
theta functions
$$\left\{\vartheta_\mu(2w)\right\}_{\mu\in
M}$$ embeds the quotient $W/\Lambda$ into projective space.
For any nonnegative integer $\ell$  let $M_{\leq \ell}$ be the
subset of
$M$ consisting of those
$\mu$ such that $\vartheta_\mu$ vanishes at the origin to order
$\leq
\ell$. We call
$$\left\{\vartheta_\mu(2w)\right\}_{\mu\in
M_{\leq 1}}$$
the {\em edited} $4\Theta$ linear system associated to the
Jacobian $(W,H,\Lambda)$. Our goal is to prove the following result.

\begin{Theorem}\label{Theorem:EditedFourTheta}
The edited
$4\Theta$ linear system embeds $W/\Lambda$ into
projective space.
\end{Theorem}
The theorem can easily be reconciled with its classical-style
formulation in the introduction  by means of the remarks
of \S\ref{subsubsection:ExplicitThetaConstruction}.
The proof of the theorem requires some preparation and is not
going to be completed until the end of
\S\ref{subsection:TangentVectorSeparation}. 

\subsubsection{Remark}
As was noted by the referee,
it is a part of the folklore of the theory of compact 
Riemann surfaces that $M_{\leq 1}=M$ generically.
Nonetheless our theorem does have some
content because in general the inequality 
$M_{\leq 1}\neq M$ does hold. 
Indeed, the latter holds ``with a vengeance'' for
hyperelliptic curves of large genus; see
\cite[Chap.\ IIIa, \S6, p.\ 3.105]{MumfordTataII}.

\subsubsection{Further remark}
We are grateful to the referee for providing us with a sketch of 
a proof of the fact that $M_{\leq 1}=M$ generically. 
We paraphrase the referee's remarks here.
On the Teichm\"{u}ller space $T$ 
classifying marked Riemann surfaces of genus $g$ 
each theta characteristic $\mu$ defines a subvariety
$V_\mu$ the points of which correspond to marked Riemann surfaces
such that $\theta_\mu$ vanishes at the origin to
order higher than parity requires.
What we have to prove
is that none of the $V_\mu$ equal $T$.
In any case, for at least one even $\mu$ and one odd $\mu$ 
we have $V_\mu\neq T$, as can be verified by the study of 
hyperelliptic curves. But since the mapping class group 
acts transitively on the possible markings of a Riemann surface,
and in turn acts transitively on the even (resp.\ odd) theta 
characteristics, we find that $V_\mu\neq T$ for all $\mu$,
which is what we wanted to prove.

\subsubsection{Recollection of the
Riemann quartic theta identity}
\label{subsubsection:QuarticRecollection} As in
Corollary~\ref{Corollary:RiemannQuarticIdentity}, let
$\{C_\mu\}_{\mu\in M}$ be the unique family of constants such that 
$$
\prod_{i=1}^4\vartheta_0\left(w_i\right)=
\sum_{\mu\in M}
C_\mu\cdot\left\{\begin{array}{cl}
&\vartheta_\mu((+w_1+w_2+w_3+w_4)/2)\\
\cdot&\vartheta_\mu((+w_1+w_2-w_3-w_4)/2)\\
\cdot&\vartheta_\mu((+w_1-w_2+w_3-w_4)/2)\\
\cdot&\vartheta_\mu((+w_1-w_2-w_3+w_4)/2)
\end{array}\right.
$$
for all $w_1,w_2,w_3,w_4\in W$. Recall that none of the $C_\mu$
vanish. 

\begin{Proposition}\label{Proposition:BasePointFree}
The 
linear system
$$\left\{\vartheta_\mu(2w)\right\}_{\mu\in
M_{\leq 0}}$$
is basepoint free (and hence so is the edited $4\Theta$
linear system).
\end{Proposition}
\proof This is well known. The method of proof in classical terms is
likely known to the reader. We sketch the proof here just in order to help
orient the reader toward our point of view. Fix
$t\in W$ and let
$w$ range over
$W$. We have
$$\frac{1}{2}\left[\begin{array}{rrrr}
1&1&1&1\\
1&1&-1&-1\\
1&-1&1&-1\\
1&-1&-1&1
\end{array}\right]
\left[\begin{array}{c}
t+w\\
t-w\\
t+w\\
t-w
\end{array}\right]=
\left[\begin{array}{c}
2t\\
0\\
2w\\
0
\end{array}\right]$$
and hence
$$\vartheta_0(t+w)^2\vartheta_0(t-w)^2\\\\ 
=\sum_{\mu\in M_{\leq 0}}C_\mu
\vartheta_\mu(2t)\vartheta_\mu(0)\vartheta_\mu(2w)
\vartheta_\mu(0).
$$
The left side of the
identity above does not vanish identically as a function of $w$ and hence
$\vartheta_\mu(2t)\neq 0$ for some
$\mu\in M_{\leq 0}$.
\qed

\subsection{The generating function $Z_t$}
We now introduce a technical device by means of which we can take maximum
advantage of the Riemann quartic theta identity.

\subsubsection{Definition}
To each $t\in W$ we associate a {\em generating function}
$$Z_t:\tilde{X}^4\times W\rightarrow \CC$$
by the rule
$$Z_t(\tilde{P}_0,\tilde{P}_1,\tilde{P}_2,\tilde{P}_3,w)
:=\left\{\begin{array}{cl}
&\vartheta_0(t-z(\tilde{P}_0)-z(\tilde{P}_1)+w)\\
\cdot&\vartheta_0(t+z(\tilde{P}_1)+z(\tilde{P}_3)-w)\\
\cdot &\vartheta_0(t-z(\tilde{P}_2)-z(\tilde{P}_3)+w)\\
\cdot & \vartheta_0(t+z(\tilde{P}_0)+z(\tilde{P}_2)-w).
\end{array}\right.
$$
The function $Z_t$ is holomorphic on $\tilde{X}^4\times W$
and also depends holomorphically on the parameter $t$.
Since 
$$\vartheta_0(-w)=\pm \vartheta_0(w),\;\;\;\vartheta_t(w)=\exp(-\pi
H(t,w))\vartheta_0(t+w),$$
we can rewrite the definition of $Z_t$ in the form
$$Z_t(\tilde{P_0},\tilde{P_1},\tilde{P_2},\tilde{P_3},w)
=\left\{\begin{array}{cl}
&\vartheta_{+t}(w-z(\tilde{P}_0)-z(\tilde{P}_1))\\
\cdot&\vartheta_{-t}(w-z(\tilde{P}_1)-z(\tilde{P}_3))\\
\cdot &\vartheta_{+t}(w-z(\tilde{P}_2)-z(\tilde{P}_3))\\
\cdot & \vartheta_{-t}(w-z(\tilde{P}_0)-z(\tilde{P}_2)).
\end{array}\right.
$$
The latter presentation of $Z_t$ is going to be quite convenient later.

\subsubsection{Key properties}
The generating function $Z_t$ viewed as a function of $t$ is a
theta function on
$W$ relative to $\Lambda$. More precisely, 
we have\\
$$\begin{array}{cl}
&Z_{t+\lambda}
(\tilde{P}_0,\tilde{P}_1,\tilde{P}_2,\tilde{P}_3,w)\\\\
=&\exp\left(4\pi H(\lambda,t)+2\pi H(\lambda,\lambda)\right)
Z_t(\tilde{P}_0,\tilde{P}_1,\tilde{P}_2,\tilde{P}_3,w)
\end{array}
$$\\
for all $\lambda\in \Lambda$. In particular, up to a nonzero constant
factor,
$Z_t$ depends only on $t\bmod{\Lambda}$.
Further and crucially, we have
$$\begin{array}{cl}
&\frac{1}{2}\left[\begin{array}{rrrr}
1&1&1&1\\
1&1&-1&-1\\
1&-1&1&-1\\
1&-1&-1&1
\end{array}\right]
\left[\begin{array}{l}
t-z(\tilde{P}_0)-z(\tilde{P}_1)+w\\
t+z(\tilde{P}_1)+z(\tilde{P}_3)-w\\
t-z(\tilde{P}_2)-z(\tilde{P}_3)+w\\
t+z(\tilde{P}_0)+z(\tilde{P}_2)-w
\end{array}\right]\\\\
&
=\left[\begin{array}{c}
2t\\\\
-z(\tilde{P}_0)+z(\tilde{P}_3)\\\\
\displaystyle 2w-\sum_{i=0}^3
z(\tilde{P}_i)\\\\
-z(\tilde{P}_1)+z(\tilde{P}_2)
\end{array}\right]
\end{array}
$$
and
hence
$$Z_t(\tilde{P}_0,\tilde{P}_1,\tilde{P}_2,\tilde{P}_3,w)=
\sum_{\mu\in M_{\leq 1}}
C_\mu \vartheta_\mu(2t)\cdot
\left\{\begin{array}{cl}
&\vartheta_\mu(-z(\tilde{P}_0)+z(\tilde{P}_3))\\\\
\cdot&\displaystyle\vartheta_\mu\left(2w
-\sum_{i=0}^3z(\tilde{P}_i)\right)\\\\
\cdot&\vartheta_\mu(-z(\tilde{P}_1)+z(\tilde{P}_2))
\end{array}\right.
$$\\
by the Riemann quartic theta identity as recollected above combined with
Corollary~\ref{Corollary:PrimeTrick}.
\subsubsection{Rationale for the terminology}
\label{subsubsection:Rationale}
We claim
that the family 
$$\left\{\begin{array}{c}
\;\\
\;\\
\;\\
\;\\
\;\\\end{array}
(\tilde{P}_0,\tilde{P}_1,\tilde{P}_2,\tilde{P}_3,w)\mapsto
\left\{\begin{array}{cl}
&\vartheta_\mu(-z(\tilde{P}_0)+z(\tilde{P}_3))\\\\
\cdot&\displaystyle\vartheta_\mu\left(2w
-\sum_{i=0}^3z(\tilde{P}_i)\right)\\\\
\cdot&\vartheta_\mu(-z(\tilde{P}_1)+z(\tilde{P}_2))
\end{array}\right.\right\}_{\mu\in M_{\leq 1}}
$$
of holomorphic functions on
$\tilde{X}^4\times W$ is $\CC$-linearly independent. 
At any rate, by Corollary~\ref{Corollary:PrimeTrick}
there exist points $\tilde{P}_0,\tilde{P}_1$
such that 
$$
\prod_{\mu\in M_{\leq 1}}\vartheta_\mu(-z(\tilde{P}_0)+z(\tilde{P}_1))\neq 0,$$
and hence by evaluation
at 
$(\tilde{P}_0,\tilde{P}_0,\tilde{P}_1,\tilde{P}_1)\in\tilde{X}^4$
we can specialize the family in question to a family of functions on
$W$ of the form 
$$\left\{D_\mu\vartheta_\mu(2(w+t_0))\right\}_{\mu\in M_{\leq 1}}$$
where the factors $D_\mu$ are nonzero constants and 
$t_0$ is some particular point of $W$. Families of the latter sort
are clearly $\CC$-linearly independent.
Thus the claim is proved.
The claim granted, it follows that for any 
$t\in W$,
to know the generating function $Z_t$ up to a nonzero constant
factor  is to know the image of
$t$ under the map to projective space determined by the edited
$4\Theta$  linear system, and vice versa. In this sense $Z_t$
packages all the information in the edited $4\Theta$ embedding
and hence deserves to be called a generating function.

\subsection{Separation of points and tangent vectors}
\label{subsection:TangentVectorSeparation}

\subsubsection{An {\em ad hoc} notion of general position}
\label{subsubsection:GeneralPosition}
Suppose we are given an integer
$n\geq g$ and a divisor $D$ of $X$ of degree $g-1$. 
We say that points $P_0,\dots,P_{n+1}\in X$ are
in {\em $(n,D)$-general position} if
$P_0,\dots,P_{n+1},\infty$ are distinct
and 
$$h^0\left(
D+n\infty-\sum_{i\in I}P_i\right)=0=h^0
\left(K-D+n\infty-\sum_{i\in I}P_i\right)$$
for all subsets $I\subset \{0,\dots,n+1\}$ of cardinality $n$,
where as usual $K$ denotes a canonical divisor of $X$. Now fix $t\in
W$ such that $\psi_D=\psi_t$, noting that $\psi_{K-D}=\psi_{-t}$ by
Corollary~\ref{Corollary:Conjugation}. Also fix distinct points
$P_0,\dots,P_{n+1}\in X\setminus\{\infty\}$ and corresponding liftings
$\tilde{P}_0,\dots,\tilde{P}_{n+1}\in
\tilde{X}$. By Corollary~\ref{Corollary:WRVT}
the points $P_0,\dots,P_{n+1}$ are in $(n,D)$-general position if and only
if
$$\prod_{\begin{subarray}{c}
I\subset\{0,\dots,n+1\}\\
\card I=n
\end{subarray}}\left(\vartheta_t\left(\sum_{i\in I}z(\tilde{P}_i)\right)
\vartheta_{-t}\left(\sum_{i\in I}z(\tilde{P}_i)\right)\right)
\neq 0.$$
It follows that the set consisting of $(n+2)$-tuples  of
points of
$X$ in
$(n,D)$-general position is open dense in $X^{n+2}$.

\begin{Proposition}\label{Proposition:SeparatesPoints} For all
$t,t'\in W$, if we have
$Z_{t}\propto Z_{t'}$, then we have $t\equiv
t'\bmod{\Lambda}$.
\end{Proposition}
\proof
Fix divisors $D$ and $D'$ of degree $g-1$ such that
$\psi_t=\psi_D$ and $\psi_{t'}=\psi_{D'}$.
Fix an integer $n\geq g+1$.
Fix points 
$P_0,\dots,P_{n+1}\in X\setminus\infty$
in both $(n,D)$-general and $(n,D')$-general position.
Fix corresponding liftings
$\tilde{P}_0,\dots,\tilde{P}_{n+1}\in
\tilde{X}$. For all indices $i,j=1,\dots,n$
we define a function $F_{ij}$ meromorphic on $\tilde{X}$ and regular on
$\tilde{X}\setminus\Gamma\tilde{\infty}$ by specializing the generating
function
$Z_t$ and multiplying by a  factor independent of $t$, as follows:
$$F_{ij}(\tilde{P}):=
\left\{\begin{array}{cl}
&\displaystyle Z_t\left(\tilde{P},\tilde{P}_i,\tilde{P}_j,
\tilde{P}_{n+1},z(\tilde{P})+\sum_{\alpha=1}^{n+1}z(\tilde{P}_\alpha)\right)\\\\
\cdot &\displaystyle\prod_{\alpha\in \{1,\dots,n\}\setminus\{i\}}
E(\tilde{P},\tilde{P}_\alpha)\cdot
\prod_{\beta\in \{1,\dots,n\}\setminus\{j\}}
E(\tilde{P},\tilde{P}_\beta).
\end{array}\right.
$$
By the second version of the definition of the generating function $Z_t$ we
have
$$F_{ij}(\tilde{P})=\left\{\begin{array}{cl}
&\displaystyle\vartheta_{+t}\left(\sum_{\alpha\in
\{1,\dots,n+1\}\setminus\{i\}}z(\tilde{P}_\alpha)\right)\\\\
\cdot&\displaystyle
\vartheta_{-t}\left(z(\tilde{P})+\sum_{\alpha\in
\{1,\dots,n\}\setminus\{i\}}z(\tilde{P}_\alpha)\right)\cdot\prod_{\alpha\in
\{1,\dots,n\}\setminus\{i\}}E(\tilde{P},\tilde{P}_\alpha)\\\\
\cdot&\displaystyle \vartheta_{+t}\left(z(\tilde{P})+\sum_{\beta\in
\{1,\dots,n\}\setminus\{j\}}z(\tilde{P}_\beta)\right)\cdot\prod_{\beta\in
\{1,\dots,n\}\setminus\{j\}}E(\tilde{P},\tilde{P}_\beta)\\\\
\cdot&\displaystyle\vartheta_{-t}\left(\sum_{\beta\in
\{1,\dots,n+1\}\setminus\{j\}}z(\tilde{P}_\beta)\right).
\end{array}\right.
$$\\
The latter presentation of $F_{ij}$ makes it clear that $F_{ij}$ does not
vanish identically.  Indeed, since the points
$P_0,\dots,P_{n+1}$ are in $(n,D)$-general position, we have
$$\prod_{i=1}^n\prod_{j=1}^n F_{ij}(\tilde{P}_0)\neq 0.$$
It is easily verified that
$F_{ij}$ is a theta function on
$\tilde{X}$ relative to
$\Gamma$ and hence represents some divisor on $X$, say $D_{ij}$. By
Proposition~\ref{Proposition:Frame} the divisor
$$D^*:=(n-1)\infty+\min_{i=1}^n D_{i1}$$ belongs to the divisor class of
$D+\infty$. Since $Z_t\propto
Z_{t'}$, a
repetition of the preceding argument proves that 
$D^*$ also belongs to the divisor class of
$D'+\infty$.   Therefore the divisors
$D$ and
$D'$ are linearly equivalent, hence
$\psi_t=\psi_D=\psi_{D'}=\psi_{t'}$ 
and hence $t\equiv t'\bmod{\Lambda}$.
\qed
\begin{Proposition}\label{Proposition:SeparatesTangentVectors}
For all $t,w_0\in W$, if the meromorphic function
$$\left.\frac{\partial}{\partial s}\log Z_{t+sw_0}\right|_{s=0}$$
on $\tilde{X}^4\times W$ reduces to a constant, then
we have $w_0=0$.
\end{Proposition}
\proof 
Put
$$
\dot{\vartheta}_{\pm t}(w):=\left.\frac{\partial}{\partial s}\log
\vartheta_{\pm t}(w+sw_0)\right|_{s=0},$$
thereby defining a meromorphic function on $W$ regular away from
the zero locus of $\vartheta_{\pm t}$. 
Note that
$$
\vartheta_{\pm(t+sw_0)}(w)=
\exp(-\pi H(\pm w_0,w)\bar{s}+\pi H(t,w_0)s)\vartheta_{\pm t}(w\pm sw_0)$$
and hence
$$\left.\frac{\partial}{\partial s}\log \vartheta_{\pm
(t+sw_0)}\right|_{s=0}=
\pi H(t,w_0)
\pm \dot{\vartheta}_{\pm t}.$$
Now fix a divisor $D$ of $X$ of degree $g-1$ such that
$\psi_t=\psi_D$,  an integer $n\geq g+1$, points
$P_0,\dots,P_{n+1}\in X\setminus\infty$
in $(n,D)$-general position and corresponding liftings
$\tilde{P}_0,\dots,\tilde{P}_{n+1}\in
\tilde{X}$.  For $i,j=1,\dots,n$ we have
$$\begin{array}{rcl}
\mbox{constant}&=&\displaystyle\left.\frac{\partial}{\partial
s}\log Z_{t+sw_0}\left(\tilde{P}_0,\tilde{P}_i,
\tilde{P}_j,\tilde{P}_{n+1},z+\sum_{\alpha=1}^{n+1}z(\tilde{P}_\alpha)
\right)\right|_{s=0}\\\\
&=&4\pi H(t,w_0)+\left\{\begin{array}{cl}
&\displaystyle\dot{\vartheta}_{+t}\left(\sum_{\alpha\in
\{1,\dots,n+1\}\setminus\{i\}}z(\tilde{P}_\alpha)\right)\\\\
-&\displaystyle
\dot{\vartheta}_{-t}\left(z+\sum_{\alpha\in
\{1,\dots,n\}\setminus\{i\}}z(\tilde{P}_\alpha)\right)\\\\
+&\displaystyle
\dot{\vartheta}_{+t}\left(z+\sum_{\beta\in
\{1,\dots,n\}\setminus\{j\}}z(\tilde{P}_\beta)\right)\\\\
-&\displaystyle\dot{\vartheta}_{-t}\left(\sum_{\beta\in
\{1,\dots,n+1\}\setminus\{j\}}z(\tilde{P}_\beta)\right).
\end{array}\right.
\end{array}
$$
We conclude that $w_0=0$ by
Proposition~\ref{Proposition:DeformationFormula}. \qed
\subsubsection{Completion of the  proof of
Theorem~\ref{Theorem:EditedFourTheta}}
By Proposition~\ref{Proposition:BasePointFree} we know at least
that the edited $4\Theta$ linear system defines a mapping from
$W/\Lambda$ to projective space.  But then, on account of
the generating function interpretation of
$ Z_t$ provided in \S\ref{subsubsection:Rationale},
we know that this mapping 
separates points by
Proposition~\ref{Proposition:SeparatesPoints} 
and 
separates tangent vectors by
Proposition~\ref{Proposition:SeparatesTangentVectors}.   \qed

\subsection{``Tying together algebraic and analytic Jacobians''}

The turn of phrase here is  borrowed from
\cite{MumfordTataII}. 

\subsubsection{Identifications}
Continuing in the setting of the proof of
Theorem~\ref{Theorem:EditedFourTheta}, we now also contemplate the
situation discussed in the first couple of paragraphs of the
introduction. Recall that in the setting of the introduction
we make the yet stronger assumption that
$$n\geq g+2.$$
Write 
$$\psi_0=\psi_{D_0}$$ for some half-canonical divisor $D_0$.
Write 
$$\EE=\OO_X(D_0+S+n\infty),\;\;\;\TT=\OO_X(T)$$ where $S$
and $T$ are divisors of
degree zero. Select
$s,t\in W$ such that
$$\psi_s=\psi_{D_0+S},\;\;\;\psi_{D_0+T}=\psi_t.$$
Make the evident identifications
$$H^0\left(X,\TT^{\pm 1}\otimes\EE\right)=\left(\begin{array}{l}
\mbox{space of theta functions }\\
\mbox{of type $(X,\infty,\psi_{s\pm t},n)$}
\end{array}\right),
$$
thereby inducing an identification of each matrix entry
$\abel(\TT)_{ij}$ with a function on $\tilde{X}^{\{0,\dots,n+1\}}$. 
We are going to write down a formula for
$\abel(\TT)_{ij}(\tilde{P}_0,\dots,\tilde{P}_{n+1})$ in terms
of the generating function $ Z_t$
and some other factors independent of $t$.
Since knowledge of $ Z_t$ up to a nonzero constant multiple
and knowledge of the image of $t$ under the map to projective space defined
by the edited $4\Theta$ linear system are equivalent, 
the formula we are going to write down ``is''
the promised factorization of the map $\TT\mapsto \abel(\TT)$
through the edited $4\Theta$ embedding.

\subsubsection{Plugging into the determinant identity}
As in the introduction, let $u$ (resp., $v$) be a row vector of length $n$
with entries forming a basis of 
$H^0(X,\TT^{-1}\otimes \EE)$ (resp., $H^0(X,\TT\otimes \EE)$)
over the ground field, the latter now taken to be $\CC$.
 Taking into account the identifications
made above, and after adjusting $u$ and $v$ by suitably
chosen nonzero constant factors, we have
$$\left(\det_{i,j=1}^n u_j^{(i)}\right)
\left(\tilde{P_1},\dots,\tilde{P_n}\right)=
\left\{\begin{array}{ccl}
&&\displaystyle\vartheta_{-t}\left(s+\sum_{i=1}^n
z(\tilde{P}_i)\right)\\\\
&\cdot&
\displaystyle\exp\left(-\pi
H\left(s,\sum_{i=1}^n
z(\tilde{P}_i)\right)\right)\\\\
&\cdot&\displaystyle\prod_{1\leq i<j\leq
n}E(\tilde{P}_i,\tilde{P}_j)
\end{array}\right.
$$\\
and
$$\left(\det_{i,j=1}^n v_j^{(i)}\right)
\left(\tilde{P_1},\dots,\tilde{P_n}\right)=
\left\{\begin{array}{ccl}
&&\displaystyle\vartheta_{t}\left(s+\sum_{i=1}^n
z(\tilde{P}_i)\right)\\\\
&\cdot&
\displaystyle\exp\left(-\pi
H\left(s,\sum_{i=1}^n
z(\tilde{P}_i)\right)\right)\\\\
&\cdot&\displaystyle\prod_{1\leq i<j\leq
n}E(\tilde{P}_i,\tilde{P}_j)
\end{array}\right.
$$\\
by plugging into the determinant identity
of Proposition~\ref{Proposition:BasicDeterminantIdentity}, it being noted
that
$$\exp(-\pi H(s,w))\vartheta_{\pm t}(s+w)$$ is a theta function
of type $(W,H,\Lambda,\psi_{s\pm t})$.
\subsubsection{Completion of the calculation}
Finally, we have\\
$$\begin{array}{ccl}
&&\abel(\TT)_{ij}
(\tilde{P}_0,\dots,\tilde{P}_{n+1})\\\\
=&&\displaystyle Z_t\left(\tilde{P_0},\tilde{P}_i,
\tilde{P}_j,\tilde{P}_{n+1},
s+\sum_{\alpha=0}^{n+1}z(\tilde{P}_\alpha)\right)
\\\\
&\cdot&\displaystyle\exp\left(-\pi
H\left(s,-2(z(\tilde{P}_0)+z(\tilde{P}_i)+z(\tilde{P}_j)
+z(\tilde{P}_{n+1}))+4\sum_{\alpha=0}^{n+1}
z(\tilde{P}_\alpha)\right)\right)\\\\
&\cdot&\displaystyle\prod_{\begin{subarray}{c}
\alpha,\beta\in \{0,\dots,n+1\}\setminus\{0,i\}\\
\alpha<\beta
\end{subarray}}E(\tilde{P}_\alpha,\tilde{P}_\beta)\cdot\prod_{\begin{subarray}{c}
\alpha,\beta\in \{0,\dots,n+1\}\setminus\{i,n+1\}\\
\alpha<\beta
\end{subarray}}E(\tilde{P}_\alpha,\tilde{P}_\beta)\\\\
&\cdot&\displaystyle\prod_{\begin{subarray}{c}
\alpha,\beta\in \{0,\dots,n+1\}\setminus\{j,n+1\}\\
\alpha<\beta
\end{subarray}}E(\tilde{P}_\alpha,\tilde{P}_\beta)
\cdot\prod_{\begin{subarray}{c}
\alpha,\beta\in \{0,\dots,n+1\}\setminus\{0,j\}\\
\alpha<\beta
\end{subarray}}E(\tilde{P}_\alpha,\tilde{P}_\beta).
\end{array}
$$\\
Note that among the factors on the right side of the identity,
only the very
first depends on
$t$. Note also that since
$n\geq g+2$, the formula above determines the generating function $ Z_t$
uniquely in terms of, say, the matrix entry $\abel(\TT)_{12}$.
The formula above generalizes 
the genus one identity
\cite[(49)]{Anderson} to arbitrary genus.


\begin{thebibliography}{99}

\bibitem[Anderson 2002]{Anderson}
Anderson, G.\ W.: {\em Abeliants and their
application to an elementary construction of
Jacobians},
Adv.\ in Math.\ \textbf{172}(2002), 169-205.

\bibitem[Fay]{Fay}
Fay, J.: {\em Theta functions on Riemann surfaces.}
Springer Lecture Notes {\bf 352}, Berlin New
York 1973.

\bibitem[Griffiths-Harris]{GriffithsHarris}
Griffith, P.; Harris, J.: {\em Principles of Algebraic Geometry.}
John Wiley and Sons, New York 1978

\bibitem[Mumford Tata I]{MumfordTataI}
Mumford, D.: {\em Tata Lectures on Theta I.}
Progr.\ Math.\ {\bf 28}(1983), Birkh\"{a}user Boston
 
\bibitem[Mumford Tata II]{MumfordTataII}
Mumford, D.: {\em Tata Lectures on Theta II.}
Progr.\ Math.\ {\bf 43}(1984), Birkh\"{a}user Boston

\bibitem[Mumford Tata III]{MumfordTataIII}
Mumford, D.: {\em Tata Lectures on Theta III.}
Progr.\ Math.\ {\bf 97}(1991), Birkh\"{a}user Boston

\bibitem[Segal-Wilson 1985]{SegalWilson}
Segal, G.; Wilson, G.: {\em Loop groups and equations of KdV type.}
Publ.\ IHES \textbf{61}(1985)5-65.

\bibitem[Weil VK]{Weil} Weil, A.: {\em Introduction
\`{a} l'\'{e}tude des vari\'{e}t\'{e}s
k\"{a}hl\'{e}riennes.} Hermann, Paris 1971


\end{thebibliography}
\end{document}